%% file: arXiv.tex
\documentclass{article}

\usepackage[T1]{fontenc}

\RequirePackage{geometry}
\usepackage[colorlinks,citecolor=blue, linkcolor = blue,urlcolor=blue]{hyperref}

\usepackage{amsthm,amsmath,amsfonts,amssymb,mathtools} 
\usepackage{bera} 
\usepackage{authblk} 
\usepackage{enumitem}
\usepackage{bbm} 
\usepackage{mathrsfs} 
\usepackage{nicefrac}
\usepackage[numbers]{natbib}

\newcommand{\iu}{{i\mkern1mu}}
\newcommand{\Rd}[0]{\mathbb{R}^{d}}

\newtheorem{theorem}{Theorem}[section]
\newtheorem{lemma}[theorem]{Lemma}
\newtheorem{corollary}[theorem]{Corollary}

\newtheorem*{theorem*}{Theorem}

\newtheorem{remark}{Remark}[section]
\newtheorem{example}{Example}[section]

\newtheorem*{acks*}{Acknowledgments}

\title{Equivalent Gaussian distributions on commutative spaces: An RKHS analysis}

\author{Michael Hediger\thanks{E-mail: \texttt{michael.hediger@math.uzh.ch}}\\ University of Zurich, Switzerland}

\begin{document}

\maketitle


\begin{abstract}
The investigation of equivalent Gaussian distributions for stochastic processes is a central problem in probability and statistics. In this context, the choice of the index set and the correlation structure, particularly their interaction, plays a crucial role. The purpose of this paper is to show how an explicit description of the corresponding reproducing kernel Hilbert space (RKHS) helps to better understand this interplay. In the stationary setting, when the index set is taken to be a homogeneous space, we show how an RKHS approach allows us to bridge the gap to harmonic analysis on commutative spaces, thereby further complementing the characterization of equivalent Gaussian distributions via their spectral measures.
\end{abstract}

\vspace{1em}
{\emergencystretch=1.5em \noindent
 \par
\vspace{0.5em}
 \par}

\input{main.tex}

\bibliographystyle{amsplain-nodash}
\bibliography{references.bib}

\begin{acks*}
I thank Serge Cohen for the encouraging and helpful discussions. I am also grateful to Luke Conners and Pavel Mnev for valuable exchanges regarding Riemannian manifolds. Further, I wish to thank Alexander Gorodnik for his advice on harmonic analysis on commutative spaces. Finally, I greatly value Reinhard Furrer's feedback and continuous support.
\end{acks*}

\end{document}

%% file: main.tex
\section{Introduction} \label{se:introduction}

Let $P_1$ and $P_2$ be two choices of Gaussian distributions for a real-valued stochastic process $X=(X_{t})_{t \in T}$. A valuable consideration is that $P_{1}$ and $P_{2}$ are either equivalent or orthogonal on $\sigma_{T}(X)$, the $\sigma$-field generated by $X$. This dichotomy has been verified by various authors (cf.\ \cite{ParzenTimeSeries}, p.\ 478). Feldman \cite{Feldman1958} provides a characterization in terms of the Gaussian spaces $H_{1}(X)$ and $H_{2}(X)$ associated with $P_1$ and $P_2$, respectively. In particular, he shows that for $P_{1}$ and $P_{2}$ to be equivalent on $\sigma_{T}(X)$ it is necessary and sufficient that there exists a linear homeomorphism $U$ from $H_{1}(X)$ onto $H_{2}(X)$ s.t.\ $U^{*}U-I$ is Hilbert--Schmidt. He refers to $U$ as an equivalence operator from $H_{1}(X)$ onto $H_{2}(X)$. Later, Rozanov \cite{Rozanov} provides an alternative proof using the entropy of the measure $P_{1}$ w.r.t.\ $P_{2}$. An approach which was pioneered earlier by H\'{a}jek~\cite{HajekJDiv, HajekOriginal}. Let $T_{n} = \{t_{1}, \dotsc, t_{n}\}$ be a finite collection of coordinates from $T$ and denote by $\sigma(Y_{n})$ the $\sigma$-field generated by the Gaussian vector $Y_{n} = (X_{t_{1}}, \dotsc, X_{t_{n}})$. Then, if $P_{\ell}^{n}$, $\ell = 1,2$, denotes the restriction of $P_{\ell}$ to $\sigma(Y_{n})$, the J-divergence (or just divergence) $J(n)$ between the finite-dimensional distributions $P_{1}^{n}$ and $P_{2}^{n}$, according to H\'{a}jek~\cite{HajekJDiv} (cf.\ Jeffreys~\cite{Jeffreys} p.\ 158, and Kullback and Leibler~\cite{KullbackLeibler}), is defined by
\begin{equation*} \label{eq:JDiv}
    J(n) = \begin{cases}E_{2}\Big[\operatorname{log}\frac{dP_{2}^{n}}{dP_{1}^{n}}\Big] - E_{1}\Big[\operatorname{log}\frac{dP_{2}^{n}}{dP_{1}^{n}}\Big], & \text{if $P_{1}^{n} \equiv P_{2}^{n} \text{ on } \sigma(Y_{n})$,} \\
    \infty, & \text{otherwise.}\end{cases}
\end{equation*}
Notice that if the covariance function $R_{\ell}$ of $X$ under $P_{\ell}$ is strictly positive definite, then $\nicefrac{dP_{2}^{n}}{dP_{1}^{n}} = \nicefrac{p_{2}^{n}(Y_{n})}{p_{1}^{n}(Y_{n})}$, with $p_{\ell}^{n}$ the $n$-dimensional normal density of $P_{\ell}^{n}$. In~\cite{HajekOriginal} H\'{a}jek shows that $P_{1}$ and $P_{2}$ are equivalent on $\sigma_{T}(X)$ if and only if $\sup_{T_{n}}J(n) < \infty$, the supremum extending over all finite subsets $T_{n}$ of $T$. A particularly simple case is when $R_{1}$ and $R_{2}$ are strictly positive definite and s.t.\ $R_{1} = \alpha^{2} R_{2}$, $\alpha \neq 1$ (e.g., $\alpha = \nicefrac{\sigma_{1}}{\sigma_{2}}$ with $X$ a Brownian motion with covariance function $\sigma_{\ell}^{2}\min\{s,t\}$ under $P_{\ell}$). In this case H\'{a}jek's result shows that $P_{1}$ and $P_{2}$ must be orthogonal on $\sigma_{T}(X)$, as can be seen by evaluating
\begin{equation*} \label{eq:JDivEval}
        2J(n) = \operatorname{tr}\big[R_{1}(n)R_{2}(n)^{-1}\big] + \operatorname{tr}\big[R_{2}(n)R_{1}(n)^{-1}\big] - 2n, 
\end{equation*}
to
\begin{equation*} \label{eq:ClosedFormJDiv}
        J(n) = \frac{1}{2}\bigg(\alpha - \frac{1}{\alpha}\bigg)^{2}n.
\end{equation*}  
In the example of the two scaled Wiener measures mentioned above, the corresponding conclusion can be traced back to the work of Cameron and Martin \cite{CameronMartin} (compare also with Example~1 in H\'{a}jek's paper \cite{HajekOriginal}). Alongside the Gaussian spaces $H_{\ell}(X)$, the RKHS $H_{T}(R_{\ell})$ for $R_{\ell}$ is generated by functions of the type $\sum_{k}\alpha_{k}R_{\ell}(\cdot, t_{k})$, $t_{k} \in T$. Characterizing equivalent Gaussian distributions via RKHS analysis is collectively referred to as the RKHS approach. Early work on the RKHS approach was done by Parzen \cite{ParzenInferenceTimeSeriesRKHSMethods,ParzenTimeSeries,Parzen1962,ParzenRKHS}, and Kallianpur and Oodaira \cite{KallianpurOodaira1} (cf.\ Oodaira \cite{Oodaira}). Later, general results were obtained by Neveu \cite{Neveu}. A comprehensive overview is given by Chatterji and Mandrekar \cite{Chatterji_Mandrekar}. For instance, for $T$ countably infinite and $X$ a zero-mean process (under $P_{1}$ and $P_{2}$), Parzen \cite{ParzenRKHS} relies on H\'{a}jek's characterization of equivalence and shows that for $P_{1}$ and $P_{2}$ to be equivalent on $\sigma_{T}(X)$ it is necessary and sufficient that the difference in covariance $R_{2} - R_{1}$ belongs to the RKHS $H_{T \times T}(R_{1} \otimes R_{2})$, defined on $T \times T$ with reproducing kernel $R_{1} \otimes R_{2}\big((s_{1},s_{2}),(t_{1},t_{2})\big) = R_{1}(s_{1},t_{1})R_{2}(s_{2},t_{2})$. Furthermore, an expression of the Radon--Nikodym derivative of $P_{2}$ w.r.t.\ $P_{1}$ is derived (cf.\ Capon \cite{Capon}). For separable $H_{T}(R_{\ell})$, Neveu relies on the inherent structure of the Gaussian spaces $H_{\ell}(X)$ and proves the following general result.
\begin{theorem*}[Neveu \cite{Neveu}, cf.\ Proposition~8.6]
    Suppose that $X$ has zero-mean function under $P_{1}$ and $P_{2}$ and the corresponding RKHSs $H_{T}(R_{1})$ and $H_{T}(R_{2})$ are separable. Then $P_{1}$ and $P_{2}$ are equivalent on $\sigma_{T}(X)$ if and only if $R_{2} - R_{1} \in H_{T \times T}(R_{1}^{2\otimes})$, the RKHS on $T\times T$ with kernel $R_{1}^{2\otimes}\big((s_{1},s_{2}),(t_{1},t_{2})\big) = R_{1}(s_{1},t_{1})R_{1}(s_{2},t_{2})$, and the eigenvalues of the Hilbert--Schmidt operator $L_{\mathcal{D}} \colon H_{T}(R_{1}) \to H_{T}(R_{1})$ associated with $\mathcal{D} = R_{2} - R_{1}$ are strictly larger than $-1$.\footnote{The operator $L_{\mathcal{D}}$ is defined pointwise via $(L_{\mathcal{D}}f)(t) = \langle L_{\mathcal{D}}f, R_{1}(\cdot, t)\rangle_{R_{1}} = \langle \mathcal{D}(\cdot, t), f\rangle_{R_{1}}$, $t \in T$, and $\langle \cdot, \cdot \rangle_{R_{1}}$ denotes the inner product on $H_{T}(R_{1})$.}      
\end{theorem*}

In the zero-mean and stationary setting, assuming $T$ to be a separable locally compact group $G$, Chow~\cite{chow1, chow2} characterizes the equivalence of $P_{1}$ and $P_{2}$ in terms of their spectral measures $\mu_{1}$ and $\mu_{2}$. In particular, using Feldman's \cite{Feldman1958} theorem, he shows that $P_{1}$ and $P_{2}$ are equivalent on $\sigma_{T}(X)$ if and only if the nonatomic parts of $\mu_{1}$ and $\mu_{2}$ agree, they share the same set of atoms $\{a_{n}\}$, and 
\begin{equation}\label{eq:ChowIntro}
    \sum_{n}d(a_{n})\bigg(1-\frac{\mu_{1}(a_{n})}{\mu_{2}(a_{n})}\bigg)^{2} < \infty,
\end{equation}
where $d(a_{n})$ is the dimension of the representation space associated with $a_{n}$. Furthermore, he extends his theorem to $G$-spaces possessing a dense orbit. In this paper, we view Chow's theorem from an RKHS perspective, thereby establishing a link to harmonic analysis on commutative spaces, where $T = G/K$ is a homogeneous space and $(G,K)$ is a Gelfand pair. This setting yields a natural complement to Chow's theorem: If $R_{\ell}$ is strictly positive definite with $K$-invariant version that is continuous and integrable, the equivalence of $P_{1}$ and $P_{2}$ on $\sigma_{T}(X)$ reduces to their agreement in the non-compact case, and to \eqref{eq:ChowIntro} in the compact case, where $\mu_{\ell}$ is replaced by the corresponding spherical transform (Theorem~\ref{thm:ChowRKHSPerspective}). Our proof relies on an explicit description of the RKHS (Lemma~\ref{lemma:RKHSCommSpaces} and Corollary~\ref{corollary:ProductRKHS}). In a supplementary analysis, we discuss the equivalence of $P_{1}$ and $P_{2}$ under potentially weaker symmetry assumptions on $T$. We follow the approach by Narcowich~\cite{Narcowich} and discuss necessary and sufficient conditions for the equivalence of $P_{1}$ and $P_{2}$ when $T$ is a compact Riemannian manifold without boundary and $R_{\ell}$ expands according to the eigenfunctions of the Laplace--Beltrami operator. Our method combines Neveu's theorem, as stated above, with Mercer's theorem (named after Mercer's work in~\cite{Mercer}). We start with a section on introductory notation and terminology, followed by a preliminary analysis, where we establish some supporting results and motivate the RKHS approach by means of two examples. 

\section{Introductory notation and terminology} \label{se:Notation}

\subsection{Primary notation}
Given a topological space $(E, \tau)$, the Borel $\sigma$-field over $E$ is denoted by $\mathfrak{B}(E)$. The space of continuous functions $f$ on $E$ is written as $C(E)$. Also, $f \in C_{c}(E)$ if $f \in C(E)$ and $f$ has compact support. For $x,y \in \Rd$, $\langle x, y \rangle = x^{\mathrm{t}}y$ identifies the dot product on $\Rd$. Given a measure space $(E,\mathcal{A},\mu)$, the space of measurable functions $f \colon E \to \mathbb{C}$ that are square-integrable on $E$ w.r.t.\ $\mu$ is denoted by $L^{2}(E,\mathcal{A},\mu)$. The canonical norm of $\langle f , g\rangle_{\mu} = \int_{E}f(x)\overline{g(x)}\mu(dx)$ on $L^{2}(E,\mathcal{A},\mu)$ is written as $\lVert \cdot \rVert_{\mu}$. If clear from the context, we make use of the short notations $L^{2}(\mu)$ or $L^{2}(E)$. The space of absolutely integrable functions (on $E$ w.r.t.\ $\mu$) is identified with $L^{1}(E)$ (resp.\ $L^{1}(\mu)$). Moreover, we write $\mu(x)$ for $\mu$ evaluated at a measurable singleton $\{x\}$. Given two measures $\mu_{1}$ and $\mu_{2}$ on $\mathcal{A}$ and $L^{2}(\mu_{1}) \supset L \subset L^{2}(\mu_{2})$, we use the notation $\lVert \cdot \rVert_{\mu_{1}} \asymp \lVert \cdot \rVert_{\mu_{2}}$ on $L$ to indicate that the norms $\lVert \cdot \rVert_{\mu_{1}}$ and $\lVert \cdot \rVert_{\mu_{2}}$ are equivalent on $L$. That is, there exist constants $\alpha_{1}, \alpha_{2} > 0$, s.t.\ for any $\varphi \in L$, $0 < \alpha_{1} \lVert \varphi \rVert_{\mu_{2}} \leq \lVert \varphi \rVert_{\mu_{1}} \leq \alpha_{2} \lVert \varphi \rVert_{\mu_{2}} < \infty$. The measures $\mu_{1}$ and $\mu_{2}$ are termed equivalent on $\mathcal{A}$ if they are mutually absolutely continuous on $\mathcal{A}$, i.e., $\mu_{1}(A) = 0$ implies $\mu_{2}(A) = 0$, $A \in \mathcal{A}$, and vice versa. If $\mu_{1}$ and $\mu_{2}$ are equivalent on $\mathcal{A}$ we write $\mu_{1} \equiv \mu_{2}$ on $\mathcal{A}$. On the other hand, $\mu_{1}$ and $\mu_{2}$ are referred to as orthogonal on $\mathcal{A}$, written as $\mu_{1} \perp \mu_{2}$ on $\mathcal{A}$, if there exists a separating set $A \in \mathcal{A}$ for which $\mu_{1}(A) = 0$ and $\mu_{2}(E \setminus A) = 0$. Generally, given a Hilbert space $\mathcal{H}$, if no particular convention is made, $\langle \cdot, \cdot \rangle_{\mathcal{H}}$ and $\lVert \cdot \rVert_{\mathcal{H}}$ will be taken as proxies for the inner product and norm on $\mathcal{H}$, respectively. Also, $\mathcal{H}^{*}$ is notation for the dual of $\mathcal{H}$. The set of strictly positive integers is denoted by $\mathbb{N}$. In the absence of ambiguity, we write $(a_{n}) = (a_{n})_{n \in \mathbb{N}}$ for a sequence indexed over $\mathbb{N}$, and $\{a_{n}\} = \{a_{n} \colon n \in \mathbb{N}\}$ for a countably infinite set. Throughout, $\ell$ defines an index set that takes the value one or two.

\subsection{Symmetric nonnegative-definite kernels and their RKHS}

Let $T$ be a set. A function $R \colon T \times T \to \mathbb{R}$ is referred to as a symmetric nonnegative-definite kernel if $R(s,t) = R(t,s)$ for all $s,t \in T$, and if for any $n \in \mathbb{N}$ and $t_{1},\dotsc,t_{n} \in T$,
\begin{equation} \label{eq:kernel2}
    \sum_{i=1}^{n}\sum_{j=1}^{n}a_{i}a_{j}R(t_{i},t_{j}) \geq 0, \quad a_{1},\dotsc,a_{n} \in \mathbb{R}.
\end{equation}
If equality in \eqref{eq:kernel2} holds only for $a_{1} = \cdots =a_{n} = 0$, $R$ is said to be strictly positive definite. Given a symmetric nonnegative-definite kernel $R$, we write $H_{T}(R)$ for the Hilbert space of real-valued functions on $T$ which satisfies
\begin{equation*} \label{eq:RKHSProp1}
    R(\cdot, t) \in H_{T}(R), \quad t \in T,
\end{equation*}
and
\begin{equation}\label{eq:RKHSRepProp}
    f(t) = \langle f, R(\cdot, t) \rangle_{R}, \quad t \in T, \quad f \in H_{T}(R), 
\end{equation}
where $\langle \cdot, \cdot \rangle_{R}$ denotes the inner product on $H_{T}(R)$. This identifies $H_{T}(R)$ as the unique RKHS with reproducing kernel $R$. The existence and uniqueness statement regarding $H_{T}(R)$ is known as the Moore--Aronszajn theorem (cf.\ \cite{Aronszajn,RKHSBook}). Property \eqref{eq:RKHSRepProp} is called the reproducing property. For $n \in \mathbb{N}$ and $T_{n} = \{t_{1}, \dotsc,t_{n}\} \subset T$, we write $R(n)$ for the $n\times n$ matrix with entries $R(t_{i},t_{j})$. In the following, it is assumed that any RKHS $H_{T}(R)$ is separable. Recall that if $T$ is a topological space and $R$ is continuous on $T \times T$, the separability of $H_{T}(R)$ is equivalent to the separability of $T$. In particular, if nothing else is mentioned, any topological space $T$ is assumed to be separable.  

\begin{example}\label{example:RKHSofMatrix}
    Let $R$ be a strictly positive-definite kernel. Then the RKHS for $R(n)$ on $T_{n}$ is identified with $\mathbb{R}^{n}$, equipped with the inner product
\begin{equation*} \label{eq:EuclideanRKHS}
    \langle v, w \rangle_{R(n)} = v^{\mathrm{t}}R(n)^{-1}w, \quad v,w \in \mathbb{R}^{n}. 
\end{equation*} 
\end{example}

\begin{example} \label{example:RKHSSquareInt}
    Let $(E,\mathcal{A},\mu)$ be a measure space with $\sigma$-finite measure $\mu$. Assume that $\gamma(t,\cdot) \in L^{2}(E,\mathcal{A},\mu)$, $t \in T$. Denote by $L_{T}(\gamma)$ the closed subspace of $L^{2}(E,\mathcal{A},\mu)$ spanned by $\{\gamma(t,\cdot) \colon t \in T\}$. Then, if $R$ is a symmetric nonnegative-definite kernel s.t.\ 
    \begin{gather*}
        R(s,t) = \int_{E}\!\gamma(s,u)\overline{\gamma(t,u)}\mu(du), \quad s,t \in T,
    \end{gather*}
    its RKHS consists of real-valued functions $f(t) = \int_{E}\!\xi_{f}(u)\gamma(t,u)\mu(du)$, $\xi_{f} \in L_{T}(\gamma)$, with inner product $\langle f, g \rangle_{R} = \langle \xi_{f} , \xi_{g} \rangle_{\mu}$.
\end{example}

It is helpful to recall Mercer's theorem and the notion of integral operators associated with Mercer kernels. For $T$ a subset of $\Rd$, a good overview is given in \cite{RKHSBook}. Regarding more general $T$, we refer to \cite{Saburou_Yoshihiro_RKHS}. If $T$ is a compact Hausdorff space equipped with a Radon measure $\mu$, we call $R \colon T \times T \to \mathbb{R}$ a Mercer kernel if it is a continuous symmetric nonnegative-definite kernel.\footnote{Since the RKHS of a Mercer kernel is separable, the separability assumption on $T$ is not required.} The integral operator $I_{R} \colon L^{2}(\mu) \to L^{2}(\mu)$ associated with a Mercer kernel $R$ is the compact and self-adjoint operator defined by  
\begin{gather*}
    [I_{R}f](t) = \int_{T}\!R(s, t)f(s)\mu(ds).
\end{gather*}
Mercer's theorem states that any Mercer kernel $R$ with corresponding integral operator $I_{R}$ attains a Mercer expansion $R(s,t) = \sum_{n}\lambda_{n}e_{n}(s)e_{n}(t)$, where $\{e_{n}\}$ is a set of orthonormal eigenfunctions of $I_{R}$ with corresponding nonzero eigenvalues $\{\lambda_{n}\}$. This has consequences on the appearance of $H_{T}(R)$. Specifically, if $R$ is a Mercer kernel with $\{e_{n}\}$ and $\{\lambda_{n}\}$ given by its Mercer expansion, then $\{\sqrt{\lambda_{n}}e_{n}\}$ is an orthonormal basis of $H_{T}(R)$, the latter consisting precisely of the functions $\sum_{n}\alpha_{n}(f)\sqrt{\lambda_{n}}e_{n}$, with inner product $\langle f,g\rangle_{R} = \sum_{n}\alpha_{n}(f)\alpha_{n}(g)$ (cf.\ Theorem~11.18 in \cite{RKHSBook} or Theorem~2.31 in \cite{Saburou_Yoshihiro_RKHS}).  

The Hilbert--Schmidt operator in Neveu's theorem is derived from the tensor product of RKHSs, which are in isometric correspondence with the corresponding RKHS of functions on $T \times T$. We review some of the basic terminology. Let $\mathcal{H}_{1}$ and $\mathcal{H}_{2}$ be two separable Hilbert spaces. Their tensor product $\mathcal{H}_{1} \otimes \mathcal{H}_{2}$, forms a new Hilbert space, generated by elementary tensors $f \otimes g$, $f \in \mathcal{H}_{1}$, $g \in \mathcal{H}_{2}$. We use the notation $\mathcal{H}_{\ell}^{2 \otimes} = \mathcal{H}_{\ell} \otimes \mathcal{H}_{\ell}$ for the two-fold tensor product of $\mathcal{H}_{\ell}$. Recall that for any element $U$ in $\mathcal{H}_{1} \otimes \mathcal{H}_{2}$, we can find a corresponding operator $A_{U} \colon \mathcal{H}_{1} \to \mathcal{H}_{2}^{*}$ s.t.\
\begin{gather*}
    (A_{U}f)(g) = \langle f \otimes g, U \rangle_{\mathcal{H}_{1} \otimes \mathcal{H}_{2}}, \quad f \in \mathcal{H}_{1},\ g \in \mathcal{H}_{2},
\end{gather*}
and $\sum_{n}\lVert A_{U}\xi_{n}\rVert^{2}_{\mathcal{H}_{2}^{*}} = \lVert U \rVert^{2}_{\mathcal{H}_{1} \otimes \mathcal{H}_{2}}$ for any orthonormal basis $\{\xi_{n}\}$ of $\mathcal{H}_{1}$. The operator $A_{U}$ is known as the Hilbert--Schmidt operator associated with $U$ (see for instance Neveu \cite{Neveu}, Proposition~6.16). By the Riesz representation theorem, we can always express
\begin{gather*}
    (A_{U}f)(g) = \langle R(A_{U}f), g \rangle_{\mathcal{H}_{2}},
\end{gather*}
where $R(A_{U}f)$ is the unique Riesz representation of the linear functional $A_{U}f$. Therefore, also the operator $L_{U} \colon \mathcal{H}_{1} \to \mathcal{H}_{2}$ that maps any $f \in \mathcal{H}_{1}$ to its Riesz representation, $L_{U}f = R(A_{U}f)$ satisfies
\begin{gather*}
    \sum_{n}\lVert L_{U}\xi_{n}\rVert^{2}_{\mathcal{H}_{2}} = \sum_{n}\lVert A_{U}\xi_{n}\rVert^{2}_{\mathcal{H}_{2}^{*}} = \lVert U \rVert^{2}_{\mathcal{H}_{1} \otimes \mathcal{H}_{2}},
\end{gather*}
where $\{\xi_{n}\}$ is any orthonormal basis $\{\xi_{n}\}$ of $\mathcal{H}_{1}$. In particular, $L_{U}$ is Hilbert--Schmidt. 

In view of two RKHS choices $H_{T}(R)$ and $H_{T}(R^{\prime})$, we denote by $H_{T \times T}(R \otimes R^{\prime})$ the (product) RKHS on $T \times T$ with kernel $R \otimes R^{\prime}\big((s_{1},s_{2}),(t_{1},t_{2})\big) = R(s_{1},t_{1})R^{\prime}(s_{2},t_{2})$. Then the two Hilbert spaces $H_{T \times T}(R \otimes R^{\prime})$ and $H_{T}(R) \otimes H_{T}(R^{\prime})$ are isometrically isomorphic. This is summarized by writing
\begin{equation} \label{eq:IsomorphicProductRKHS}
    \widetilde{H_{T}(R) \otimes H_{T}(R^{\prime})} = H_{T \times T}(R \otimes R^{\prime}).
\end{equation}  
A proof is given by Aronszajn \cite{Aronszajn} (cf.\ Section~8). We also refer to Section~5.5 in \cite{RKHSBook}. In fact, $H_{T \times T}(R \otimes R^{\prime})$ consists precisely of the functions $\widetilde{U}(s,t) = \sum_{i}\sum_{j}\alpha_{i,j}e_{i}(s)e_{j}^{\prime}(t)$ that satisfy $\sum_{i,j}\alpha_{i,j}^{2} < \infty$. In the latter, $\{e_{i}\}$ and $\{e_{j}^{\prime}\}$ are orthonormal bases in $H_{T}(R)$ and $H_{T}(R^{\prime})$, respectively, and the sum of squared coefficients equals the squared norm of $\widetilde{U}$. Using the correspondence \eqref{eq:IsomorphicProductRKHS}, we can make the identification $L_{U} = L_{\widetilde{U}}$ for the Hilbert--Schmidt operator associated with an element $\widetilde{U} \in H_{T \times T}(R \otimes R^{\prime})$. Specifically, for the two-fold tensor product choice $H_{T}(R)^{2 \otimes}$, $(A_{U}f)(g) = \langle L_{\widetilde{U}}f, g \rangle_{R}$ and by the reproducing property of $H_{T}(R)$,
\begin{align*}
    (L_{\widetilde{U}}f)(t) &= \langle L_{\widetilde{U}}f, R(\cdot, t) \rangle_{R} = (A_{U}f)(R(\cdot, t)) \\
    &= \langle f \otimes R(\cdot, t), U \rangle_{H_{T}(R)^{2 \otimes}} = \langle fR(\cdot, t), \widetilde{U} \rangle_{R^{2\otimes}}.
\end{align*}
In particular, if $f = R(\cdot, s)$, i.e., $fR(\cdot, t) = R^{2\otimes}(\cdot, (s,t))$, we obtain $(L_{\widetilde{U}}R(\cdot, s))(t) = \widetilde{U}(s,t)$. This makes,
\begin{gather*}
    (L_{\widetilde{U}}R(\cdot, s))(t) = \langle L_{\widetilde{U}}R(\cdot, s), R(\cdot, t) \rangle_{R} = \widetilde{U}(s,t) = \langle \widetilde{U}(\cdot,t), R(\cdot, s) \rangle_{R},
\end{gather*}
and therefore, by the separability of $H_{T}(R)$, 
\begin{equation} \label{eq:DefHSOPRKHS}
    (L_{\widetilde{U}}f)(t) = \langle L_{\widetilde{U}}f, R(\cdot, t) \rangle_{R} = \langle \widetilde{U}(\cdot,t),f \rangle_{R}, \quad f \in H_{T}(R).
\end{equation}
From here, it is readily seen that $L_{\widetilde{U}} \colon H_{T}(R) \to H_{T}(R)$ is self-adjoint, i.e., $\langle L_{\widetilde{U}}f, g \rangle_{R} = \langle f, L_{\widetilde{U}}g \rangle_{R}$ for any $f,g \in H_{T}(R)$, if and only if $\widetilde{U}(s,t) = \widetilde{U}(t,s)$.

\subsection{Conventions for the stochastic process}
For the remainder of this article, it is assumed that $X$ has zero-mean function under both measures $P_{1}$ and $P_{2}$. That is, $P_{1}$ and $P_{2}$ differ only in terms of their covariance functions $R_{1}$ and $R_{2}$. We also refer to $P_{\ell}$ as a centered Gaussian distribution on $\sigma_{T}(X)$ with covariance function $R_{\ell}$. The linear span of $\{X_{t} \colon t \in T\}$ is denoted by $H_{0}(X)$, i.e., the Gaussian space $H_{\ell}(X)$ identifies its closure in $L^{2}(P_{\ell})$. Finally, we reserve the notation $\mathcal{D}$ for the difference in covariance $R_{2} - R_{1}$ defined on $T \times T$.

\section{Preliminary analysis} \label{se:PreliminaryResults}

The following theorem is attributed to Aronszajn and Neveu, it is an interpretation of Neveu's theorem from the perspective of Aronszajn's work \cite{Aronszajn} (the statement is formulated in the article by Chatterji and Mandrekar \cite{Chatterji_Mandrekar}).  

\begin{theorem}[Aronszajn--Neveu]\label{thm:Aronszajn-Neveu}
    The zero-mean Gaussian distributions $P_{1}$ and $P_{2}$ are equivalent on $\sigma_{T}(X)$ if and only if $H_{T}(R_{1}) = H_{T}(R_{2})$ and $\mathcal{D} \in H_{T \times T}(R_{1}^{2\otimes})$.
\end{theorem}

For completeness, we state and verify two lemmas that streamline the proof of the theorem.

\begin{lemma} \label{lemma:EQUAL_RKHS_EQUAL_GAUSSSPACE}
    In order that $\lVert \cdot \rVert_{P_{1}} \asymp \lVert \cdot \rVert_{P_{2}}$ on $H_{0}(X)$ it is necessary and sufficient that $H_{T}(R_{1}) = H_{T}(R_{2})$.
\end{lemma}

\begin{proof}
    This is a consequence of Aronszajn’s inclusion theorem (cf.\ Corollary~IV$_{3}$ on p.\ 383 in his work~\cite{Aronszajn}). In particular, $H_{T}(R_{1}) = H_{T}(R_{2})$ if and only if there exist constants $\alpha_{1},\alpha_{2} > 0$ s.t.\ the matrices $\alpha_{2}^{2}R_{2}(n)-R_{1}(n)$ and $\alpha_{1}^{2}R_{1}(n)-R_{2}(n)$ are positive-semidefinite. Therefore, if $a= (a_{1}, \dotsc, a_{n}) \in \mathbb{R}^{n}$ and $t_{1}, \dotsc, t_{n}$ is a finite collection of coordinates from $T$, we make use of the correspondence   
\begin{gather*}
    a^{\mathrm{t}}R_{\ell}(n)a = \big\lVert \sum_{i=1}^{n}a_{i}R_{\ell}(\cdot,t_{i})\big \rVert_{R_{\ell}}^{2} = \big\lVert \sum_{i=1}^{n}a_{i}X_{t_{i}} \big \rVert_{P_{\ell}}^{2}, \quad \ell = 1,2,
\end{gather*}
and the statement of the lemma is verified. 
\end{proof}

Certainly, a violation of $\lVert \cdot \rVert_{P_{1}} \asymp \lVert \cdot \rVert_{P_{2}}$ on $H_{0}(X)$ leads to orthogonal measures $P_{1}$ and $P_{2}$ on $\sigma_{T}(X)$, so that the implication $H_{T}(R_{1}) \neq  H_{T}(R_{2})$ $\Rightarrow$ $P_{1}\perp P_{2}$ on $\sigma_{T}(X)$ can be concluded from Lemma~\ref{lemma:EQUAL_RKHS_EQUAL_GAUSSSPACE}. The analogous observation, with $H_{T}(R_{\ell})$ replaced with the Cameron--Martin space of $P_{\ell}$, is given in Proposition~2.7.3 of Bogachev's book \cite{Bogachev}. As an example of two orthogonal Gaussian measures that share an identical RKHS, we can again think of the two scaled Wiener measures mentioned in the introduction.

\begin{lemma} \label{lemma:SpekrumHSAndEqualRKHS}
    Under the assumption that $\mathcal{D} \in H_{T \times T}(R_{1}^{2 \otimes})$, $H_{T}(R_{1})$ equals $H_{T}(R_{2})$ if and only if the eigenvalues of the operator $L_{\mathcal{D}}$ as defined in \eqref{eq:DefHSOPRKHS} are strictly larger than $-1$.  
\end{lemma}

\begin{proof}
    By Lemma~\ref{lemma:EQUAL_RKHS_EQUAL_GAUSSSPACE} and Neveu's theorem, it remains to show that if the two RKHS are equal, the eigenvalues of $L_{\mathcal{D}}$ are strictly larger than $-1$. Therefore, assume that $H_{T}(R_{1}) = H_{T}(R_{2})$. Let $M= I+L_{\mathcal{D}}$, where $I \colon H_{T}(R_{1}) \to H_{T}(R_{1})$ is the identity. By definition of $L_{\mathcal{D}}$, 
    \begin{gather*}
        MR_{1}(\cdot, t) = R_{1}(\cdot, t) + L_{\mathcal{D}}R_{1}(\cdot, t) = R_{1}(\cdot, t) + \mathcal{D}(\cdot, t) = R_{2}(\cdot, t). 
    \end{gather*}
    Since $H_{T}(R_{1}) = H_{T}(R_{2})$, $M$ has corresponding inverse $M^{-1} \colon H_{T}(R_{1}) \to H_{T}(R_{1})$ that satisfies $M^{-1}R_{2}(\cdot, t) = R_{1}(\cdot, t)$. Certainly $M^{-1}$ is bounded, specifically, it is self-adjoint and Hilbert--Schmidt. In particular, there is a constant $c >0$, the constant given by the equivalence of the norms $\lVert \cdot \rVert_{R_{1}}$ and $\lVert \cdot \rVert_{R_{2}}$, s.t.\ for any $f \in H_{T}(R_{1})$, $\langle Mf, f\rangle_{R_{2}} \geq c \lVert f \rVert^{2}_{R_{2}}$ (i.e., $M \geq cI$). If we choose $f = e_{i}$ an eigenfunction of $L_{\mathcal{D}}$ with corresponding eigenvalue $\lambda_{i}$, the latter is equivalent to $\lambda_{i} \geq c-1 > -1$.  
\end{proof}

Using two examples, we show below how Theorem~\ref{thm:Aronszajn-Neveu} can be applied to examine the equivalence or orthogonality of $P_{1}$ and $P_{2}$. The first gives an original take on the situation where $X$ is stationary under $P_{\ell}$ and $T$ corresponds to the entire $\Rd$ (cf.\ Chapter~III of Ibragimov's and Rozanov's book \cite{GaussianRandomProcesses}, using a spectral approach). The second complements the analysis by Arafat et al.~\cite{GaussianSphereEquivOrtho} (cf.\ Yadrenko~\cite{Yadrenko}), done on the unit sphere. In both examples, we rely on an explicit description of the RKHS.

\begin{example}[\textbf{Stationary processes on real coordinate spaces}] \label{example:StationaryRealCoordinateSpaces}
Let $T = \Rd$ and assume that
\begin{gather*}
    R_{\ell}(s+h,t+h) = R_{\ell}(s,t), \quad s,t,h \in \Rd.
\end{gather*}
That is, $X$ is stationary under $P_{\ell}$. Define $k_{\ell}(t) = R_{\ell}(t,0)$, $t \in \Rd$. We observe that $R_{\ell}(s,t) = k_{\ell}(s-t)$. Suppose that $k_{\ell}$ is continuous at zero (i.e., continuous everywhere). Then, by Bochner's theorem (for real coordinate spaces) \cite{Bochner}, 
\begin{gather*}
    R_{\ell}(s,t) = \int_{\Rd}\!\operatorname{e}^{\iu \langle s, \lambda  \rangle}\overline{\operatorname{e}^{\iu \langle t, \lambda  \rangle}}F_{\ell}(d\lambda), \quad s,t \in \Rd,
\end{gather*}
for some finite measure $F_{\ell}$, uniquely defined on $\mathfrak{B}(\Rd)$. Notice that we are in the framework of Example~\ref{example:RKHSSquareInt}. Actually, if $R_{\ell}$ is strictly positive definite and $k_{\ell} \in L^{1}(\Rd)$, it is known that $H_{\Rd}(R_{\ell})$ consists of continuous functions $f \in L^{2}(\Rd)$ s.t.\ $\nicefrac{\hat{f}}{\sqrt{\hat{k}_{\ell}}} \in L^{2}(\Rd)$ with inner product,
\begin{gather*}
    \langle f , g \rangle_{R_{\ell}} = (2\pi)^{\nicefrac{-d}{2}}\int_{\Rd}\!\frac{\hat{f}(\lambda)\overline{\hat{g}(\lambda)}}{\hat{k}_{\ell}(\lambda)}d\lambda.
\end{gather*}
In the latter, $f \mapsto \hat{f}$ denotes the Fourier transform of $f$ (cf.\ Wendland \cite{Wendland_Book}, Theorem~10.12). In particular, $H_{\Rd}(R_{\ell}) \subset L^{2}(\Rd)$. Upon the isometric correspondence given in \eqref{eq:IsomorphicProductRKHS}, it follows that $H_{\Rd \times \Rd}(R_{\ell}^{2\otimes}) \subset L^{2}(\Rd \times \Rd)$. Then, since the map $t \mapsto \int_{\Rd}\!\mathcal{D}(s,t)^{2}ds$ is constant, we conclude that $\mathcal{D}$ does not belong to $H_{\Rd \times \Rd}(R_{1}^{2\otimes})$, unless $R_{1}$ and $R_{2}$ are equal. This shows that for $T = \Rd$, the Gaussian measures $P_{1}$ and $P_{2}$ are orthogonal on $\sigma_{\Rd}(X)$ as soon as $R_{1} \neq R_{2}$.\footnote{For an overview of results concerning the characterization of equivalent Gaussian distributions for real stationary processes, we refer to the books by Yadrenko \cite{Yadrenko}, Ibragimov and Rozanov \cite{GaussianRandomProcesses}, and Gikhman and Skorokhod \cite{StochasticProcesses}.}
\end{example}

\begin{example}[\textbf{Isotropic processes on the sphere}] \label{example:CovSphere}

Let $T = \mathbb{S}^{d-1}$, $d \geq 3,$ be the unit sphere in $\Rd$. Denote by $\vartheta \colon \mathbb{S}^{d-1} \times \mathbb{S}^{d-1} \to [0,\pi]$ the great-circle (geodesic) distance on $\mathbb{S}^{d-1}$, i.e., $\vartheta(s,t) = \operatorname{arccos}(\langle s, t \rangle)$. Assume that 
\begin{gather*}
    R_{\ell}(s,t) = \psi_{\ell}(\vartheta(s,t)), \quad s,t \in \mathbb{S}^{d-1},
\end{gather*}
where $\psi_{\ell} \colon [0,\pi] \to \mathbb{R}$ is continuous and s.t.\ $\psi_{\ell}(0) > 0$. Given $n \in \mathbb{N} \cup \{0\}$, write $S_{n}^{i}$, $i = 1, \dotsc, h(n)$, for the spherical harmonics of degree $n$ (cf.\ Chapter~XI (Section~11.3) in \cite{Erdelyi} or also Chapter~IV (Section~2) in \cite{SteinWeiss}). Then, in the sense of Schoenberg \cite{Schoenberg}, the following series representation of $R_{\ell}$ is valid, 
\begin{equation} \label{eq:ExSphere2}
        R_{\ell}(s,t) = \sum_{n=0}^{\infty}\sum_{i=1}^{h(n)}\lambda_{n}(\ell)S_{n}^{i}(s)S_{n}^{i}(t), \quad \lambda_{n}(\ell) = \frac{a_{\ell}(n)}{h(n)}.
\end{equation}
In the latter, the coefficients $a_{\ell}(n)$ are strictly positive for infinitely many $n$ (cf.\ Theorem~1 in \cite{GneitingMultiquadric} or also Section~5 of Chapter~I in \cite{Yadrenko}). In fact, $R_{\ell}$ is a Mercer kernel and the spherical harmonics $\{S_{n}^{i}\}$ form an orthonormal basis of $L^{2}(\mathbb{S}^{d-1})$, consisting of eigenfunctions of the corresponding integral operator $I_{R_{\ell}}$ with eigenvalues $\{\lambda_{n}(\ell)\}$. Therefore, we know that $H_{\mathbb{S}^{d-1}}(R_{\ell})$ consists of functions 
\begin{equation} \label{eq:ElementsRKHSSphere}
        f(t) = \sum_{n=0}^{\infty}\sum_{i=1}^{h(n)}\beta_{n,i}(f)\sqrt{\lambda_{n}(\ell)}S_{n}^{i}(t),
\end{equation}
with inner product,
\begin{equation} \label{eq:InnerProductRKHSSphere}
    \langle f, g \rangle_{R_{\ell}} = \sum_{n=0}^{\infty}\sum_{i=1}^{h(n)}\beta_{n,i}(f)\beta_{n,i}(g), \quad f,g \in H_{\mathbb{S}^{d-1}}(R_{\ell}).
\end{equation}
Let $H_{\mathbb{S}^{d-1}}(R_{1}) = H_{\mathbb{S}^{d-1}}(R_{2})$ and $\mathcal{D} \in H_{\mathbb{S}^{d-1} \times \mathbb{S}^{d-1}}(R_{1}^{2\otimes})$. Since the two RKHSs are equal, we write $\mathcal{D} \in H_{\mathbb{S}^{d-1} \times \mathbb{S}^{d-1}}(R_{2}^{2\otimes})$ in terms of its orthonormal basis expansion
\begin{equation*} \label{eq:ProductRKHSSphere}
    \mathcal{D}(s,t) = \sum_{n=0}^{\infty}\sum_{i=1}^{h(n)}\sum_{m=0}^{\infty}\sum_{j=1}^{h(m)}\gamma_{i,j}(n,m)\sqrt{\lambda_{n}(2)\lambda_{m}(2)}S_{n}^{i}(s)S_{m}^{j}(t),
\end{equation*}
with
\begin{gather*}
    \lVert \mathcal{D} \rVert_{R_{2}^{2\otimes}}^{2} = \sum_{n,m}\sum_{i,j}\gamma_{i,j}(n,m)^{2} < \infty.
\end{gather*}
If we compare coefficients with \eqref{eq:ExSphere2}, we observe that
\begin{equation} \label{eq:CoefficientsSphere}
    \gamma_{i,j}(n,m) = 
    \begin{cases}
        \frac{\lambda_{n}(2)-\lambda_{n}(1)}{\lambda_{n}(2)} = \frac{a_{2}(n)-a_{1}(n)}{a_{2}(n)}, & \text{$n = m$ and $i=j$,} \\
            0, & \text{otherwise.} 
    \end{cases}
\end{equation}
This makes
\begin{equation} \label{eq:EquivalenceSphere}
    \lVert \mathcal{D} \rVert_{R_{2}^{2\otimes}}^{2} = \sum_{n=0}^{\infty}h(n)\bigg(1-\frac{a_{1}(n)}{a_{2}(n)}\bigg)^{2}.
\end{equation}
Hence, by Theorem~\ref{thm:Aronszajn-Neveu}, $P_{1} \equiv P_{2}$ on $\sigma_{\mathbb{S}^{d-1}}(X)$ implies that the sum on the right-hand side of Equation~\eqref{eq:EquivalenceSphere} is finite. Regarding the converse, if the series in \eqref{eq:EquivalenceSphere} converges, we define $\gamma_{i,j}(n,m)$ as in \eqref{eq:CoefficientsSphere} and conclude by the definition of the product RKHS that $\mathcal{D}\in H_{\mathbb{S}^{d-1} \times\mathbb{S}^{d-1}}(R_{2}^{2\otimes})$. Further, by the convergence of the series in \eqref{eq:EquivalenceSphere}, we must have
\begin{gather*}
    \lim_{n \to \infty}\frac{a_{1}(n)}{a_{2}(n)} = \lim_{n \to \infty}\frac{a_{2}(n)}{a_{1}(n)} = 1.
\end{gather*}
In particular, the sequences $(\nicefrac{a_{1}(n)}{a_{2}(n)})$ and $(\nicefrac{a_{2}(n)}{a_{1}(n)})$ are bounded. If $f$ belongs to $H_{\mathbb{S}^{d-1}}(R_{1})$, it follows from \eqref{eq:ElementsRKHSSphere} and \eqref{eq:InnerProductRKHSSphere} that $f = \sum_{n=0}^{\infty}\sum_{i=1}^{h(n)}\beta_{n,i}(f)\sqrt{\lambda_{n}(1)}S_{n}^{i}$ with $\sum_{n=0}^{\infty}\sum_{i=1}^{h(n)}\beta_{n,i}(f)^{2} < \infty$. Hence, 
\begin{gather*}
        f = \sum_{n=0}^{\infty}\sum_{i=1}^{h(n)}\beta_{n,i}(f)\sqrt{\frac{a_{1}(n)}{a_{2}(n)}}\sqrt{\lambda_{n}(2)}S_{n}^{i} \in H_{T}(R_{2}). 
\end{gather*}
That is, $H_{\mathbb{S}^{d-1}}(R_{1}) \subset H_{\mathbb{S}^{d-1}}(R_{2})$. By the same reasoning, the reverse inclusion holds. Therefore, $\mathcal{D} \in H_{\mathbb{S}^{d-1} \times \mathbb{S}^{d-1}}(R_{1}^{2\otimes})$. In conclusion:
\begin{gather*}
    \text{$P_{1} \equiv P_{2}$ on $\sigma_{\mathbb{S}^{d-1}}(X)$} \quad \Leftrightarrow \quad \sum_{n=0}^{\infty}h(n)\bigg(1-\frac{a_{1}(n)}{a_{2}(n)}\bigg)^{2} < \infty.
\end{gather*}
This gives a complementary RKHS approach that reproduces Theorem~1 in \cite{GaussianSphereEquivOrtho}.
\end{example}

\section{Stationary processes on commutative spaces}

Taking an algebraic perspective, the two-sphere $\mathbb{S}^{2}$ is a $G$ space under the transitive action of the rotation group $SO(3)$. In particular, we can identify $\mathbb{S}^{2}$ with the quotient space $SO(3)/K_{p} = \{gK_{p} \colon g \in SO(3)\}$, the set of left cosets of $K_{p}$ in $SO(3)$, where $K_{p}$ is the stabilizer subgroup of an arbitrary point on the sphere (the subgroup of rotations that leave $p$ untouched). The identification of the two spaces is given by the $G$-set homeomorphism $\Phi$ that maps a coset $gK_{p}$ to the corresponding point $g\cdot p = gp$ on the sphere. In summary, this is written as $\mathbb{S}^{2} \cong SO(3)/SO(2)$. As an example, let $s = g_{s}K_{p}$ and $t = g_{t}K_{p}$ be two elements of $SO(3)/K_{p}$ and define
\begin{gather*}
    R_{\ell}(s,t) = \psi_{\ell}(\vartheta(\Phi(s),\Phi(t))) = \psi_{\ell}(\vartheta(g_{s}p, g_{t}p)), 
\end{gather*}
where $\psi_{\ell}$ is given as in Example~\ref{example:CovSphere}. Then $R_{\ell}$ is an instance of an $SO(3)$-invariant kernel on $SO(3)/K_{p} \times SO(3)/K_{p}$. That is, for any $g \in SO(3)$, $R_{\ell}(gs,gt) = R_{\ell}(s,t)$. If we choose to study $T = SO(3)/K_{p}$, $X$ becomes $G$-invariant (stationary) under $P_{\ell}$. A treatment on group invariant stochastic processes is given in Malyarenko's book \cite{Malyarenko} (cf.\ the earlier work by Yaglom~\cite{YaglomInvariant} or Askey and Bingham~\cite{Askey_Bingham}). For a comprehensive algebraic and geometric overview, the reader is referred to Wolf~\cite{Wolf}, Folland~\cite{Folland}, or Helgason \cite{Helgason}.

Throughout this section, $G$ is assumed to be a locally compact group. In general, we write $g \cdot g^{\prime} = gg^{\prime}$ for the group operation and $e = gg^{-1}$ identifies the identity element on $G$. A continuous function $\varphi \colon G \to \mathbb{C}$ is said to be positive definite if for any $n \in \mathbb{N}$ and $g_{1}, \dotsc, g_{n} \in G$,
\begin{equation*} \label{eq:PDFun}
\sum_{i=1}^{n}\sum_{j=1}^{n}z_{i}\overline{z}_{j}\varphi(g_{i}g_{j}^{-1}) \geq 0, \quad z_{1}, \dotsc, z_{n} \in \mathbb{C}.
\end{equation*}
As motivation, it is useful to recall the setting were $G$ is a locally compact abelian (LCA) group, i.e., Bochner's theorem applies. 

\begin{remark}[Positive-definite functions on LCA groups]
If $G$ is LCA, it is common to use $+$ for the group operation and $0$ for the identity element. Recall that the dual group $G^{*}$ consists of continuous characters of $G$, i.e., $\chi \in G^{*}$ if $\chi \colon G \to \mathbb{C}$ is continuous, $\lvert \chi(g) \rvert = 1$, $g \in G$, and 
\begin{equation*} \label{eq:CharacterRole}
    \chi(g+g^{\prime}) = \chi(g)\chi(g^{\prime}), \quad g,g^{\prime} \in G.
\end{equation*}
In particular, $\chi(0) = 1$ and $\chi(-g) = \overline{\chi(g)}$, $g \in G$. The following is known as Bochner's theorem \cite{Rudin}: A continuous function $\varphi \colon G \to \mathbb{C}$ is positive definite if and only if,
\begin{equation} \label{eq:Bochner}
    \varphi(g) = \int_{G^{*}}\chi(g)\mu(d\chi), \quad g \in G,
\end{equation}
for some finite measure $\mu$, uniquely defined on $\mathfrak{B}(G^{*})$. If $T = G$ and $X$ is stationary under $P_{\ell}$, $R_{\ell}(s+g,t+g) = R_{\ell}(s,t)$, $s,t,g \in G$, it follows that the function $k_{\ell} \colon G \to \mathbb{C}$, defined by $k_{\ell}(t) = R_{\ell}(t,0)$, $t \in G$, is positive definite. If $k_{\ell}$ is continuous, we deduce from \eqref{eq:Bochner} that there exists a finite measure $\mu_{\ell}$ on $\mathfrak{B}(G^{*})$ s.t.\  
\begin{gather*}
    R_{\ell}(s,t) = \int_{G^{*}}\chi(s)\overline{\chi(t)}\mu_{\ell}(d\chi), \quad s,t \in G.
\end{gather*}
Having in mind Example~\ref{example:StationaryRealCoordinateSpaces}, if $G = \Rd$, the dual group $G^{*}$ is isomorphic to $\Rd$ with isomorphism $\lambda \mapsto \operatorname{exp}(\iu \langle \cdot, \lambda \rangle)$. That is, any continuous character $\chi$ of $\Rd$ is given by $\chi(\cdot) = \chi_{\lambda}(\cdot) = \operatorname{exp}(\iu \langle \cdot, \lambda \rangle)$ for some $\lambda \in \Rd$. If we identify $\Lambda \colon G^{*} \to \Rd$, $\Lambda(\chi_{\lambda}) = \lambda$, it follows from \eqref{eq:Bochner} that any continuous positive-definite function $\varphi \colon \Rd \to \mathbb{C}$ admits a representation $\varphi(x) = \int_{\Rd}\!\operatorname{exp}(\iu \langle x, \lambda \rangle)F(d \lambda)$ where $F$ is the pushforward of $\mu$ by $\Lambda$.
\end{remark}

The analysis of equivalent Gaussian distributions on groups goes back to the work of Grenander \cite{Grenander} and Chow~\cite{chow1, chow2} (cf.\ \cite{Chatterji_Mandrekar}). Chow's analysis is first carried out on locally compact groups and then extended to $G$-spaces possessing a dense orbit. His approach is based on Feldman's characterization via equivalence operators between Gaussian spaces. 

For the remainder of this section, we assume that $T = G/K$ is a homogeneous space\footnote{A transitive $G$-space that is isomorphic to $G/K$.}, where $K$ is a compact subgroup of $G$. It is assumed that $X$ is $G$-invariant (or stationary) under $P_{\ell}$, i.e., for any $g \in G$ and $s,t \in T$, $R_{\ell}(gs,gt) = R_{\ell}(s,t)$. We require $G$ to be unimodular. That is, the Haar measure is both left and right invariant w.r.t.\ the group action. Also, $\int_{K}\!dk = 1$, i.e., the Haar measure is normalized on $K$ (integration on $G$ against the Haar measure is written as $dg$). The goal is to reformulate Chow's result (see Chow's theorem below) to the setting where $T$ is a commutative space and $R_{\ell}$ is subject to classical Fourier regularity conditions. In particular, having access to the Fourier transform, we show how an RKHS approach fits naturally within the harmonic analysis setting (cf.\ the earlier Example~\ref{example:StationaryRealCoordinateSpaces}).

Recall that it is possible to recover Bochner's characterization \eqref{eq:Bochner} for $G$ not necessarily abelian. This aligns with the study of spherical functions (rooted in the work of Cartan \cite{Cartan} and Weyl \cite{Weyl}, see also the articles by Godement~\cite{GodementSpherical} and Tamagawa \cite{Tamagawa}). A complex-valued function $f$ defined on $G$ is called bi-$K$-invariant if $f(kgk^{\prime}) = f(g)$ for all $g \in G$ and $k,k^{\prime} \in K$. The set of bi-$K$-invariant members of $C_{\mathrm{c}}(G)$ is denoted by ${C}_{\mathrm{c}}(K \backslash G/K)$. The latter is regarded as an algebra over $\mathbb{C}$, with multiplication
\begin{gather*}
    f_{1} * f_{2} (v) = \int_{G}\!f_{1}(g)f_{2}(g^{-1}v)dg. 
\end{gather*}
If ${C}_{\mathrm{c}}(K \backslash G/K)$ is commutative, we call $(G,K)$ a Gelfand pair (following the work of Gelfand \cite{Gelfand}). In this case, $T$ is said to be commutative, which we assume henceforth. A function $\zeta \colon G \to \mathbb{C}$ is called spherical for $(G,K)$, if it is continuous, bi-$K$-invariant, $\zeta(e) = 1$, and for any $f \in {C}_{\mathrm{c}}(K \backslash G/K)$, $\zeta$ satisfies $f * \zeta = \lambda_{\scriptscriptstyle f} \zeta$ for some $\lambda_{\scriptscriptstyle f} \in \mathbb{C}$. Let $\mathscr{L}$ denote the space of positive-definite spherical functions for the Gelfand pair $(G,K)$. We view $\mathscr{L}$ as a topological space, with topology given by the compact-open topology. Then, if a continuous bi-$K$-invariant function $\varphi \colon G \to \mathbb{C}$ is positive definite, there exists a finite measure $\mu$, uniquely defined on $\mathfrak{B}(\mathscr{L})$, s.t.\
\begin{equation} \label{eq:GelfandGodementBochner1}
    \varphi(g) = \int_{\mathscr{L}}\!\zeta(g)\mu(d\zeta).
\end{equation}
The latter representation is due to Godement \cite{Godement1958} (cf.\ the earlier work of Gelfand \cite{Gelfand}). By the $G$-invariance of $R_{\ell}$, it follows that $R_{\ell}(s,t) = R_{\ell}(g_{t}^{-1}g_{s}K, K)$, $s,t \in T$, $s = g_{s}K$, $t = g_{t}K$. Therefore, if we define $\varphi_{\ell}(g) = R_{\ell}(gK,K)$, $\varphi_{\ell}$ is bi-$K$-invariant and positive definite. We refer to $\varphi_{\ell}$ as the $K$-invariant version of $R_{\ell}$. If $\varphi_{\ell}$ happens to be continuous, we deduce from \eqref{eq:GelfandGodementBochner1} that
\begin{equation} \label{eq:InvariantKernel1}
    R_{\ell}(s,t) = \varphi_{\ell}(g_{t}^{-1}g_{s}) = \int_{\mathscr{L}}\!\zeta(g_{t}^{-1}g_{s})\mu_{\ell}(d\zeta), \quad s= g_{s}K, \ t =g_{t}K.
\end{equation}

\begin{remark} \label{remark:InvariantKernel1}
    Given the choices for $G$ and $K$, besides the requirements for $R_{\ell}$ (resp. $\varphi_{\ell}$), the only assumption underlying \eqref{eq:InvariantKernel1} is that $(G,K)$ is a Gelfand pair. Specifically, if $G$ is a Lie group, any (Riemannian) symmetric pair $(G,K)$ is a Gelfand pair. This result is due to Gelfand \cite{Gelfand} (cf.\ \cite{Helgason}, p.\ 408). More generally, it is known that any weakly symmetric space (cf.\ Selberg~\cite{Selberg}) gives rise to a Gelfand pair.
\end{remark}

\begin{remark} \label{remark:InvariantKernel2}
    Reconsidering our introductory setting, where $G = SO(3)$ and $K = K_{p}$ is the stabilizer subgroup of an arbitrary point on $\mathbb{S}^{2}$, the kernel of Example~\ref{example:CovSphere} reads as 
    \begin{gather*}
        \psi_{\ell}(\vartheta(s^{\prime},t^{\prime})) = \int_{\mathscr{L}}\!\zeta(g_{t^{\prime}}^{-1}g_{s^{\prime}})\mu_{\ell}(d\zeta), \quad s^{\prime},t^{\prime} \in \mathbb{S}^{2}, \ s^{\prime}= g_{s^{\prime}}p, \ t^{\prime} =g_{t^{\prime}}p.
    \end{gather*}
    The connection between the spherical harmonics series expansion \eqref{eq:ExSphere2} and the latter integral is provided by compactness (see for instance \cite{Helgason}).
\end{remark}

Given $\zeta \in \mathscr{L}$, it is always possible to find a corresponding (spherical) irreducible unitary representation $\pi_{\zeta}$ on a Hilbert space $H_{\zeta}$ with cyclic vector $u_{\zeta}$. That is, for any $g \in G$, $\zeta(g) = \langle \pi_{\zeta}(g)u_{\zeta}, u_{\zeta} \rangle_{\zeta}$, where $\langle \cdot , \cdot \rangle_{\zeta}$ denotes the inner product on $H_{\zeta}$ (see for instance Theorem~3.4, Chapter~IV, \S 3, in \cite{Helgason}). If $H_{\zeta}$ is finite dimensional, $d(\zeta)$ is the dimension of $H_{\zeta}$, otherwise $d(\zeta) = \infty$. The following is due to Chow.

\begin{theorem*}[Chow \cite{chow2}] \label{thm:Chow}
    Assume that $\varphi_{\ell}$ is continuous, i.e., \eqref{eq:InvariantKernel1} holds. Then the Gaussian distributions $P_{1}$ and $P_{2}$ are equivalent on $\sigma_{T}(X)$ if and only if
\begin{enumerate}[label=(\roman*)] 
\item \label{item:Chow1} the nonatomic parts of $\mu_{1}$ and $\mu_{2}$ agree;
\item \label{item:Chow2} $\mu_{1}$ and $\mu_{2}$ share the same set of atoms $\{a_{n}\}$ and $\sum_{n}d(a_{n})\Big(1-\frac{\mu_{1}(a_{n})}{\mu_{2}(a_{n})}\Big)^{2} < \infty$.
\end{enumerate}
\end{theorem*}

To view Chow's theorem from an RKHS perspective, we first want to make the RKHS accessible. To begin, it is helpful to recall the notion of the spherical transform. Let $BS$ be the subset of spherical functions for $(G,K)$ that are bounded. The space of functions $f$ that belong to $L^{1}(G)$ and are bi-$K$-invariant is denoted by $L^{1}(K \backslash G/K)$. Then the spherical transform is the map $\mathcal{S} \colon L^{1}(K \backslash G/K) \to BS$ defined by, 
\begin{gather*}
    [\mathcal{S}(f)](\zeta) = \int_{G}\!f(g)\zeta(g^{-1})dg, \quad \zeta \in BS.
\end{gather*} 
Its inversion formula is linked to \eqref{eq:GelfandGodementBochner1} via the Plancherel measure $\mu_{\scriptscriptstyle P}$ on $\mathfrak{B}(\mathscr{L})$ (cf.\ Theorem~9.4.1 in~\cite{Wolf}). Specifically, let $B(K \backslash G/K)$ denote the Bochner space for $G$, all linear combinations of positive-definite and bi-$K$-invariant members of $C(G)$, then $f \in B(K \backslash G/K) \cap L^{1}(K \backslash G/K)$ implies that
\begin{gather*}
    f(g) = \int_{\mathscr{L}}\![\mathcal{S}(f)](\zeta)\zeta(g)\mu_{\scriptscriptstyle P}(d\zeta), \quad g \in G.
\end{gather*}
Having in mind $\varphi_{\ell}$, which belongs to $B(K \backslash G/K) \cap L^{1}(K \backslash G/K)$ if it is continuous and integrable, we can see from the proof of Theorem~9.4.1 in~\cite{Wolf} that up to a set of $\mu_{\scriptscriptstyle P}$ measure zero,
\begin{equation} \label{eq:SphericalTransformPlancherelRDDerivative}
    [\mathcal{S}(\varphi_{\ell})](\zeta) = \frac{\mu_{\ell}(d\zeta)}{\mu_{\scriptscriptstyle P}(d \zeta)}.
\end{equation}
In particular, given the required assumptions on $\varphi_{\ell}$, $\mu_{\scriptscriptstyle P}$ a.e., $\mathcal{S}(\varphi_{\ell})$ is real-valued and nonnegative. Also, for $R_{\ell}$ strictly positive definite, $\mathcal{S}(\varphi_{\ell})$ can not vanish on a set of $\mu_{\scriptscriptstyle P}$ positive measure. We will see in the following how the expression \eqref{eq:SphericalTransformPlancherelRDDerivative} is a key tool for our analysis on $T$.

In terms of notation, we use a Fraktur font $\mathfrak{f}$ for the pullback $g \mapsto f(gK)$ of a function $f$ on $T$. Also, for any $t \in T$, $g_{t}$ denotes the group element s.t.\ $t = g_{t}K$. The Bochner space associated with $T$ is denoted by $B(T)$, i.e., $f \in B(T)$ if $\mathfrak{f}$ can be written as a linear combination of positive-definite and right-$K$-invariant members of $C(G)$. Similarly, we use the notation $f \in L^{1}(T)$ (resp.\ $f \in L^{2}(T)$) if $\mathfrak{f} \in L^{1}(G)$ (resp.\ $\mathfrak{f} \in L^{2}(G)$). Accordingly, we identify $\lVert f \rVert_{L^{1}(T)} = \lVert \mathfrak{f} \rVert_{L^{1}(G)}$ and $\lVert f \rVert_{L^{2}(T)} = \lVert \mathfrak{f} \rVert_{L^{2}(G)}$. Then, if $f \in L^{1}(T)$, we have access to its Fourier transform (cf.\ Definition 9.6.7 in~\cite{Wolf}) 
\begin{equation} \label{eq:FourierTransformDef}
    [\mathcal{F}(f)](\zeta) = \dot{\pi}_{\zeta}(\mathfrak{f})u_{\zeta}, \quad \dot{\pi}_{\zeta}(\mathfrak{f})u_{\zeta} = \int_{G}\!\mathfrak{f}(g)[\pi_{\zeta}(g)u_{\zeta}]dg.
\end{equation}
In particular, if $f \in B(T) \cap L^{1}(T)$, Fourier inversion holds (cf.\ Proposition~9.6.5 in \cite{Wolf})
\begin{gather*}
    f(t) = \int_{\mathscr{L}}\!\big\langle [\mathcal{F}(f)](\zeta), \pi_{\zeta}(g_{t})u_{\zeta} \big\rangle_{\zeta}\mu_{\scriptscriptstyle P}(d\zeta).
\end{gather*}
In addition, Plancherel's theorem tells us that, viewed as a function on $L^{1}(T) \cap L^{2}(T)$, the Fourier transform $\mathcal{F}$ extends to an isometry ${}^{\scriptscriptstyle \uparrow \!}\mathcal{F}$ of $L^{2}(T)$ onto the corresponding $L^{2}$ direct integral Hilbert space $\mathcal{H}^{2}(G,K)$ derived from $(\mathscr{L}, \mathfrak{B}(\mathscr{L}), \mu_{\scriptscriptstyle P})$, $\{H_{\zeta} \colon \zeta \in\mathscr{L}\}$, and the sections $\dot{\pi}_{\zeta}(\mathfrak{f})u_{\zeta}$, $\mathfrak{f} \in C_{c}(G)$ (cf.\ Corollary 9.6.6 and subsequent Theorem~9.6.12 in \cite{Wolf}). In particular,
\begin{gather*}
    \lVert {}^{\scriptscriptstyle \uparrow \!}\mathcal{F}(f) \rVert_{\mathcal{H}^{2}(G,K)} = \int_{\mathscr{L}}\!\lVert [{}^{\scriptscriptstyle \uparrow \!}\mathcal{F}(f)](\zeta) \rVert^{2}_{\zeta}\mu_{\scriptscriptstyle P}(d\zeta) = \lVert f \rVert^{2}_{L^{2}(T)}, \quad f \in L^{2}(T).
\end{gather*}
In terms of notation, although the extension ${}^{\scriptscriptstyle \uparrow \!}\mathcal{F}$ might not be given by the formula \eqref{eq:FourierTransformDef}, we shall not distinguish between $\mathcal{F}$ and ${}^{\scriptscriptstyle \uparrow \!}\mathcal{F}$ in the subsequent.

Recall that $R_{\ell}(s, t) = \varphi(g_{t}^{-1}g_{s})$ so that $\mathfrak{R}_{\ell}(\cdot, t) = \varphi(g_{t}^{-1}\cdot)$ and therefore, by the left-invariance of the Haar measure on $G$, $\int_{G}\mathfrak{R}_{\ell}(g, t)dg = \int_{G}\varphi_{\ell}(g_{t}^{-1}g)dg = \int_{G}\varphi_{\ell}(g)dg$. Specifically, $\varphi_{\ell} \in C(G) \cap L^{1}(G)$ if and only if $R_{\ell}(\cdot, t) \in C(T) \cap L^{1}(T)$.  

\begin{lemma}\label{lemma:RKHSCommSpaces}
    Let $R_{\ell}$ be strictly positive definite and $\varphi_{\ell} \in C(G) \cap L^{1}(G)$. Then the space $H_{T}(R_{\ell})$ consists of functions $f \in C(T) \cap L^{2}(T)$ satisfying
    \begin{equation} \label{eq:MembersOfRKHS}
       f(t) = \int_{\mathscr{L}}\!\big\langle [\mathcal{F}(f)](\zeta), \pi_{\zeta}(g_{t})u_{\zeta} \big\rangle_{\zeta}\mu_{\scriptscriptstyle P}(d\zeta),
    \end{equation}
    with inner product,
    \begin{equation} \label{eq:InnerProductRKHS}
    \langle f_{1}, f_{2} \rangle_{R_{\ell}} = \int_{\mathscr{L}}\!\frac{\big\langle [\mathcal{F}(f_{1})](\zeta), [\mathcal{F}(f_{2})](\zeta) \big\rangle_{\zeta}}{[\mathcal{S}(\varphi_{\ell})](\zeta)}\mu_{\scriptscriptstyle P}(d\zeta), \quad f_{1},f_{2} \in H_{T}(R_{\ell}).
    \end{equation}
\end{lemma}

\begin{proof}
  Let $\mathcal{H}(\ell)$ be composed of functions of type \eqref{eq:MembersOfRKHS}. Using \eqref{eq:SphericalTransformPlancherelRDDerivative}, we see that $\mathcal{H}(\ell)$ equipped with \eqref{eq:InnerProductRKHS} defines an inner product space. We show that it is complete. Let $(f_{n})$ be a Cauchy sequence w.r.t.\ the norm $\lVert \cdot \rVert_{\mathcal{H}(\ell)}$ induced by \eqref{eq:InnerProductRKHS}. Using Plancherel's theorem, the sequence $(\nicefrac{\mathcal{F}(f_{n})}{\sqrt{\mathcal{S}(\varphi_{\ell})}})$ is a Cauchy sequence in $\mathcal{H}^{2}(G,K)$. Since the latter is complete, it has a limit $a$ in $\mathcal{H}^{2}(G,K)$. Define $\xi = a\sqrt{\mathcal{S}(\varphi_{\ell})}$. By assumption, $\varphi_{\ell} \in L^{1}(G)$, so that $\mathcal{S}(\varphi_{\ell}) \in BS$. In particular, $\xi \in \mathcal{H}^{2}(G,K)$. Also, by H\"{o}lder's inequality and the fact that $\mathcal{S}(\varphi_{\ell}) \in L^{1}(\mu_{\scriptscriptstyle P})$, $\xi \in \mathcal{H}^{1}(G,K)$ (the $L^{1}$ direct integral, see for instance Definition~9.6.7 in \cite{Wolf}). Define $f = \mathcal{\mathcal{F}}^{-1}(\xi)$. By the Plancherel theorem, $f \in L^{2}(T)$. Then, since $\xi \in \mathcal{H}^{1}(G,K) \cap \mathcal{H}^{2}(G,K)$, we see that $f$ is continuous and pointwise recovered as
  \begin{gather*}
      f(t) = \int_{\mathscr{L}}\!\big\langle \xi(\zeta), \pi_{\zeta}(g_{t})u_{\zeta} \big\rangle_{\zeta}\mu_{\scriptscriptstyle P}(d\zeta) =\int_{\mathscr{L}}\!\big\langle [\mathcal{F}(f)](\zeta), \pi_{\zeta}(g_{t})u_{\zeta} \big\rangle_{\zeta}\mu_{\scriptscriptstyle P}(d\zeta).
  \end{gather*}
  Moreover, we have that
  \begin{gather*}
      \lVert f_{n} - f \rVert^{2}_{\mathcal{H}(\ell)} = \bigg\lVert \frac{\mathcal{F}(f_{n})}{\sqrt{\mathcal{S}(\varphi_{\ell})}} - a \bigg\rVert^{2}_{\mathcal{H}^{2}(G,K)} \xrightarrow[]{ n \to \infty } 0,
  \end{gather*}
  which shows that $\mathcal{H}_{\ell}$ equipped with \eqref{eq:InnerProductRKHS} is complete. As noted above, since $\varphi_{\ell}$ is by assumption continuous and absolutely integrable, the same holds for $R_{\ell}(\cdot,t)$ for any $t \in T$. Also, for any $k \in K$ and $g \in G$, $\mathfrak{R}_{\ell}(gk, t) = \mathfrak{R}_{\ell}(g, t)$, so that $\mathfrak{R}_{\ell}$ is right-$K$-invariant. Moreover, because $\varphi_{\ell}$ is positive definite, so is $\mathfrak{R}_{\ell}(\cdot, t)$. In particular, $R_{\ell}(\cdot,t)$ is a member of the Bochner space. Hence, using Fourier inversion, 
  \begin{gather*}
      R_{\ell}(s,t) = \int_{\mathscr{L}}\!\big\langle [\mathcal{F}(R_{\ell}(\cdot,t))](\zeta), \pi_{\zeta}(g_{s})u_{\zeta} \big\rangle_{\zeta}\mu_{\scriptscriptstyle P}(d\zeta).
  \end{gather*}
  Also, upon \eqref{eq:GelfandGodementBochner1}, by \eqref{eq:SphericalTransformPlancherelRDDerivative}, we see that for any $t \in T$,
  \begin{gather*}
      R_{\ell}(s,t) = \int_{\mathscr{L}}\!\zeta(g_{t}^{-1}g_{s})\mu_{\ell}(d\zeta) = \int_{\mathscr{L}}\! \big\langle [\mathcal{S}(\varphi_{\ell})](\zeta) \pi_{\zeta}(g_{t})u_{\zeta}, \pi_{\zeta}(g_{s})u_{\zeta} \big\rangle_{\zeta} \mu_{\scriptscriptstyle P}(d\zeta),
  \end{gather*}
  which identifies $\mu_{\scriptscriptstyle P}$ a.e., $[\mathcal{F}(R_{\ell}(\cdot,t))](\zeta) = [\mathcal{S}(\varphi_{\ell})](\zeta) \pi_{\zeta}(g_{t})u_{\zeta}$. Thus, 
  \begin{align*}
      \lVert R_{\ell}(\cdot,t) \rVert_{\mathcal{H}(\ell)}^{2} &= \int_{\mathscr{L}}\!\frac{\big\langle [\mathcal{S}(\varphi_{\ell})](\zeta) \pi_{\zeta}(g_{t})u_{\zeta}, [\mathcal{S}(\varphi_{\ell})](\zeta) \pi_{\zeta}(g_{t})u_{\zeta} \big\rangle_{\zeta}}{[\mathcal{S}(\varphi_{\ell})](\zeta)}\mu_{\scriptscriptstyle P}(d\zeta) \\
      &= \int_{\mathscr{L}}\![\mathcal{S}(\varphi_{\ell})](\zeta)\big\langle \pi_{\zeta}(g_{t})u_{\zeta}, \pi_{\zeta}(g_{t})u_{\zeta} \big\rangle_{\zeta}\mu_{\scriptscriptstyle P}(d\zeta) \\
      &= \int_{\mathscr{L}}\![\mathcal{S}(\varphi_{\ell})](\zeta)\big\langle \pi_{\zeta}(e)u_{\zeta}, u_{\zeta} \big\rangle_{\zeta}\mu_{\scriptscriptstyle P}(d\zeta) \\
      &= \int_{\mathscr{L}}\![\mathcal{S}(\varphi_{\ell})](\zeta) \zeta(e)\mu_{\scriptscriptstyle P}(d\zeta) = \lVert \mathcal{S}(\varphi_{\ell}) \rVert_{L^{1}(\mu_{\scriptscriptstyle P})} < \infty.
  \end{align*}
  This shows that $R_{\ell}(\cdot,t) \in \mathcal{H}(\ell)$ for any $t \in T$. In view of the identification of $\mathcal{F}(R_{\ell}(\cdot,t))$, we readily check that the inner product \eqref{eq:InnerProductRKHS} satisfies the reproducing property. In summary, by the Moore--Aronszajn theorem, $\mathcal{H}(\ell) = H_{T}(R_{\ell})$.    
\end{proof}

\begin{corollary}\label{corollary:ProductRKHS}
    Let $R_{\ell}$ be strictly positive definite and $\varphi_{\ell} \in C(G) \cap L^{1}(G)$. Then the members of $H_{T\times T}(R_{\ell}^{2\otimes})$ belong to $C(T \times T) \cap L^{2}(T \times T)$ and are representable as 
    \begin{gather*}
        f(s,t) = \int_{\mathscr{L}^{2}}\!\big\langle [\mathcal{F}_{2\otimes}(f)](\zeta_{1},\zeta_{2}), [\pi_{\zeta_{1}}(g_{s})u_{\zeta_{1}}] \otimes [\pi_{\zeta_{2}}(g_{t})u_{\zeta_{2}}]\big\rangle_{\!\scriptscriptstyle{\zeta_{1} \otimes \zeta_{2}}}\mu_{\scriptscriptstyle P}\!\otimes\!\mu_{\scriptscriptstyle P}[(d(\zeta_1, \zeta_{2})],
    \end{gather*}
    where $[\mathcal{F}_{2\otimes}(f)](\zeta_{1},\zeta_{2}) = \sum_{i,j}a_{i,j}[\mathcal{F}(e_{i})](\zeta_{1}) \otimes [\mathcal{F}(e_{j})](\zeta_{2})$, with $\sum_{i,j}\alpha_{i,j}^{2} < \infty$, for some orthonormal basis $\{e_{n}\}$ of $H_{T}(R_{\ell})$ and $\langle \cdot, \cdot \rangle_{\!\scriptscriptstyle{\zeta_{1} \otimes \zeta_{2}}}$ is the inner product on $H_{\zeta_{1}} \otimes H_{\zeta_{2}}$. In addition, the inner product on $H_{T\times T}(R_{\ell}^{2\otimes})$ satisfies
    \begin{equation} \label{eq:InnerProduct_ProductRKHS}
        \langle f_{1}, f_{2}\rangle_{R_{\ell}^{2\otimes}} = \int_{\mathscr{L}^{2}}\! \frac{\big\langle [\mathcal{F}_{2\otimes}(f_{1})](\zeta_{1},\zeta_{2}), [\mathcal{F}_{2\otimes}(f_{2})](\zeta_{1},\zeta_{2}) \big\rangle_{\!\scriptscriptstyle{\zeta_{1} \otimes \zeta_{2}}}}{[\mathcal{S}(\varphi_{\ell})](\zeta_{1})[\mathcal{S}(\varphi_{\ell})](\zeta_{2})}\mu_{\scriptscriptstyle P}\!\otimes\!\mu_{\scriptscriptstyle P}[d(\zeta_1, \zeta_{2})].
    \end{equation}
\end{corollary}

\begin{proof}
    The fact that the functions of $H_{T\times T}(R_{\ell}^{2\otimes})$ belong to $L^{2}(T \times T)$ is derived from the isometric relationship \eqref{eq:IsomorphicProductRKHS}. Also, since the kernel $R_{\ell}^{2\otimes}$ is continuous on $T^{2} \times T^{2}$, $H_{T\times T}(R_{\ell}^{2\otimes})$ consists of continuous functions only. Let $\{e_{n}\}$ be an orthonormal basis of of $H_{T}(R_{\ell})$ and write an arbitrary element $f$ of $H_{T\times T}(R_{\ell}^{2\otimes})$ as $f(s,t) = \sum_{i,j}\alpha_{i,j}e_{i}(s)e_{j}(t)$ where $\sum_{i,j}\alpha_{i,j}^{2} < \infty$. For simplicity, we identify,
    \begin{gather*}
        A_{(\zeta_{1},\zeta_{2})}(s,t) = \big\langle [\mathcal{F}(e_{i})](\zeta_{1}), [\mathcal{F}(R_{\ell}(\cdot, s))](\zeta_{1})\big\rangle_{\zeta_{1}}\big\langle [\mathcal{F}(e_{j})](\zeta_{2}), [\mathcal{F}(R_{\ell}(\cdot, t))](\zeta_{2})\big\rangle_{\zeta_{2}},
    \end{gather*}
    and $B_{(\zeta_{1},\zeta_{2})}(s,t) = [\pi_{\zeta_{1}}(g_{s})u_{\zeta_{1}}] \otimes [\pi_{\zeta_{2}}(g_{t})u_{\zeta_{2}}]$. Then, by Lemma~\ref{lemma:RKHSCommSpaces},
    \begin{align*}
        f(s,t) &= \lim_{n\to \infty}\sum_{i,j=1}^{n}\alpha_{i,j}\big\langle e_{i}e_{j}, R_{\ell}^{2\otimes}(\cdot, (s,t))\big\rangle_{R_{\ell}^{2\otimes}} \\
        &= \lim_{n\to \infty}\sum_{i,j=1}^{n}\alpha_{i,j}\big\langle e_{i}, R_{\ell}(\cdot, s)\rangle_{R_{\ell}}\langle e_{j}, R_{\ell}(\cdot, t)\big\rangle_{R_{\ell}} \\
        &=\lim_{n\to \infty}\sum_{i,j=1}^{n}\alpha_{i,j}\int_{\mathscr{L}^{2}}\!\frac{A_{(\zeta_{1},\zeta_{2})}(s,t)}{[\mathcal{S}(\varphi_{\ell})](\zeta_{1})[\mathcal{S}(\varphi_{\ell})](\zeta_{2})}\mu_{\scriptscriptstyle P}\!\otimes\!\mu_{\scriptscriptstyle P}[(d(\zeta_1, \zeta_{2})](d\zeta) \\
        &= \lim_{n\to \infty}\sum_{i,j=1}^{n}\alpha_{i,j}\int_{\mathscr{L}^{2}}\!\big\langle [\mathcal{F}(e_{i})](\zeta_{1}) \otimes [\mathcal{F}(e_{j})](\zeta_{2}), B_{(\zeta_{1},\zeta_{2})}(s,t) \big\rangle_{\!\scriptscriptstyle{\zeta_{1} \otimes \zeta_{2}}}\mu_{\scriptscriptstyle P}\!\otimes\!\mu_{\scriptscriptstyle P}[(d(\zeta_1, \zeta_{2})] \\
        &= \lim_{n\to \infty}\int_{\mathscr{L}^{2}}\! \big\langle \sum_{i,j=1}^{n}\alpha_{i,j}[\mathcal{F}(e_{i})](\zeta_{1}) \otimes [\mathcal{F}(e_{j})](\zeta_{2}), B_{(\zeta_{1},\zeta_{2})}(s,t) \big\rangle_{\!\scriptscriptstyle{\zeta_{1} \otimes \zeta_{2}}}\mu_{\scriptscriptstyle P}\!\otimes\!\mu_{\scriptscriptstyle P}[d(\zeta_1, \zeta_{2})] \\
        &= \int_{\mathscr{L}^{2}}\!\big\langle [\mathcal{F}_{2\otimes}(f)](\zeta_{1},\zeta_{2}), B_{(\zeta_{1},\zeta_{2})}(s,t)\big\rangle_{\!\scriptscriptstyle{\zeta_{1} \otimes \zeta_{2}}}\mu_{\scriptscriptstyle P}\!\otimes\!\mu_{\scriptscriptstyle P}[(d(\zeta_1, \zeta_{2})].
    \end{align*}
    Observe that interchanging the limit and the integral is justified by the fact that $\sum_{i,j}\alpha_{i,j}^{2}$ is finite and $R_{\ell}^{2\otimes}(\cdot, (s,t)) \in H_{T\times T}(R_{\ell}^{2\otimes})$. The formula \eqref{eq:InnerProduct_ProductRKHS} is derived in a similar way.  
\end{proof}

\begin{remark}\label{remark:Fourier_L1_GTimesG}
    In the special case where $f \in L^{1}(T \times T)$, we can take $\mathcal{F}_{2\otimes}$ from the previous corollary as the natural tensor product version of \eqref{eq:FourierTransformDef}
    \begin{gather*}
        [\mathcal{F}_{2\otimes}(f)](\zeta_{1},\zeta_{2}) = \int_{G \times G}\!\mathfrak{f}(g_{1},g_{2})[\pi_{\zeta_{1}}(g_{1})u_{\zeta_{1}}] \otimes [\pi_{\zeta_{2}}(g_{2})u_{\zeta_{2}}]d[(g_{1},g_{2})].
    \end{gather*}
\end{remark}

Our RKHS analysis yields the following. 

\begin{theorem}\label{thm:ChowRKHSPerspective}
    Let $R_{\ell}$ be strictly positive definite and $\varphi_{\ell} \in C(G) \cap L^{1}(G)$. Then, 
    \begin{enumerate}[label=(\roman*)] 
      \item \label{item:ChowRKHS1} if $G$ is not compact, $P_{1} \equiv P_{2}$ on $\sigma_{T}(X)$ if and only if $\varphi_{1} = \varphi_{2}$;
      \item \label{ChowRKHS2} if $G$ is compact, $\mu_{1}$ and $\mu_{2}$ have the same set of atoms $\{a_{n}\}$ and $P_{1} \equiv P_{2}$ on $\sigma_{T}(X)$ if and only if $\sum_{n}d(a_{n})\big(1-\frac{[\mathcal{S}(\varphi_{1})](a_{n})}{[\mathcal{S}(\varphi_{2})](a_{n})}\big)^{2} < \infty$.
\end{enumerate}
\end{theorem}

\begin{proof}
    The difference $\mathcal{D}(s,t)$ is given by $[\varphi_{2}-\varphi_{1}](g_{t}^{-1}g_{s})$. By Theorem~\ref{thm:Aronszajn-Neveu}, for $P_{1}$ and $P_{2}$ to be equivalent on $\sigma_{T}(X)$ it is necessary that $\mathcal{D} \in H_{T\times T}(R_{1}^{2\otimes})$. By Corollary~\ref{corollary:ProductRKHS}, the latter implies that $\mathcal{D} \in L^{2}(T \times T)$. If $G$ is not compact, this is true if and only if $\varphi_{1} = \varphi_{2}$. Hence, item~\ref{item:ChowRKHS1} is verified. For the remainder of the proof, assume that $G$ is compact. Then, for any $\zeta \in \mathscr{L}$, $d(\zeta) < \infty$, and the Plancherel measure $\mu_{\scriptscriptstyle P}$ is purely atomic (cf.\ Chapter~5 in \cite{Wolf}). By \eqref{eq:SphericalTransformPlancherelRDDerivative}, it follows that $\mu_{\ell}$ is purely atomic. Thus, using the uniqueness of $\mu_{\scriptscriptstyle P}$, $\mu_{1}$ and $\mu_{2}$ share the same set of atoms $\{a_{n}\}$. In particular, up to $\mu_{\scriptscriptstyle P}$ measure zero,
    \begin{equation*} \label{eq:FourierPlancherelSpectralRelation}
       [\mathcal{S}(\varphi_{1})](\zeta) = \frac{\mu_{1}(d \zeta)}{\mu_{2}(d \zeta)}[\mathcal{S}(\varphi_{2})](\zeta) \quad \Leftrightarrow \quad [\mathcal{S}(\varphi_{2} - \varphi_{1})](\zeta) = \bigg(1- \frac{\mu_{1}(d \zeta)}{\mu_{2}(d \zeta)}\bigg)[\mathcal{S}(\varphi_{2})](\zeta).
    \end{equation*}
    To prove item~\ref{ChowRKHS2}, we apply Theorem~\ref{thm:Aronszajn-Neveu} along with Corollary~\ref{corollary:ProductRKHS}. Assume w.l.o.g.\ that $G$ has Haar measure one. By Plancherel's theorem, using the $G$-invariance of $R_{\ell}$, we see that
    \begin{gather*}
        \int_{\mathscr{L}^{2}}\! \big\lVert [\mathcal{F}_{2\otimes}(\mathcal{D})](\zeta_1, \zeta_{2})\big\rVert_{\!\scriptscriptstyle{\zeta_{1} \otimes \zeta_{2}}}^{2}\mu_{\scriptscriptstyle P}\!\otimes\!\mu_{\scriptscriptstyle P}[d(\zeta_1, \zeta_{2})] = \lVert \mathcal{D} \rVert^{2}_{L^{2}(T \times T)} = \int_{G}\!\big(\varphi_{2}(g) - \varphi_{1}(g)\big)^{2}dg.
    \end{gather*}
    In particular, we identify,
    \begin{gather*}
        [\mathcal{F}_{2\otimes}( \mathcal{D})](\zeta_{1},\zeta_{2}) = \frac{[\mathcal{S}(\varphi_{2} - \varphi_{1})](\zeta_{1})}{\sqrt{\mu_{\scriptscriptstyle P}(\zeta_{1})}}\mathbbm{1}_{\{ \zeta_{1}\}}(\zeta_{2})\big[\pi_{\zeta_{1}}(e)u_{\zeta_{1}} \otimes \pi_{\zeta_{2}}(e)u_{\zeta_{2}}\big].
    \end{gather*}
    Thus, keeping in mind \eqref{eq:InnerProduct_ProductRKHS}, we can formally write,
    \begin{equation}\label{eq:TargetTheorem}
        \lVert \mathcal{D} \rVert_{R_{2}^{2\otimes}}^{2} = \sum_{n \in \mathbb{N}}d(a_{n})\bigg(1-\frac{[\mathcal{S}(\varphi_{1})](a_{n})}{[\mathcal{S}(\varphi_{2})](a_{n})}\bigg)^{2}.
    \end{equation}
    For simplicity, let $S(1,2)$ denote the series on the right-hand side of the previous equation. If $P_{1} \equiv P_{2}$ on $\sigma_{T}(X)$, if follows from Theorem~\ref{thm:Aronszajn-Neveu} that $\mathcal{D} \in H_{T\times T}(R_{2}^{2\otimes})$, i.e., $S(1,2)$ converges. Regarding the converse, assume that $S(1,2)$ converges. By Fourier inversion, we already know that for any $t \in T$, $R_{2}(\cdot,t)$ takes the form \eqref{eq:MembersOfRKHS}. Also, the convergence of $S(1,2)$ implies that there exists a strictly positive constant $c$ s.t.\ 
    \begin{gather*}
        \lVert R_{2}(\cdot,t) \rVert_{R_{1}}^{2} \leq c\lVert R_{2}(\cdot,t) \rVert_{R_{2}}^{2}.
    \end{gather*}
    That is, $H_{T}(R_{2}) \subset H_{T}(R_{1})$. Since the reverse inclusion is also satisfied, it follows that $H_{T}(R_{1}) = H_{T}(R_{2})$. Then, by \eqref{eq:TargetTheorem}, since $S(1,2)$ converges, $\mathcal{D} \in H_{T\times T}(R_{2}^{2\otimes})$, and by the equality of the RKHSs, $\mathcal{D} \in H_{T\times T}(R_{1}^{2\otimes})$. Hence, by Theorem~\ref{thm:Aronszajn-Neveu}, $P_{1} \equiv P_{2}$ on $\sigma_{T}(X)$.
\end{proof}

\section{On Mercer kernels for compact Riemannian manifolds without boundary} \label{se:CompactM}

Let $T$ be a compact Hausdorff space with Radon measure $\mu$ and take $R_{\ell}$ to be a Mercer kernel. For the moment, let us assume that we return to the setting of the preceding section. The symmetry of $T$, as reflected in the commutativity of the convolution algebra ${C}_{\mathrm{c}}(K \backslash G/K)$, intertwines naturally with the $G$-invariant kernel $R_{\ell}$ via its integral operator. According to Folland's~\cite{Folland} notation, using the right-$K$-invariance of $\mathfrak{R}_{\ell}(\cdot,t)$, we observe that
\begin{gather*}
    [P\mathfrak{R}_{\ell}(\cdot,t)](gK) = \int_{K}\!\mathfrak{R}_{\ell}(gk,t)dk = \mathfrak{R}_{\ell}(g,t)\int_{K}\!dk = R_{\ell}(gK, t).
\end{gather*}
So that (cf.\ \cite{Folland}, Theorem~2.51 and subsequent corollary),  
\begin{gather*}
    \int_{T}R_{\ell}(gK, t)\mu(dgK) = \int_{T}\![P\mathfrak{R}_{\ell}(\cdot,t)](gK)\mu(dgK) = \int_{G}\varphi_{\ell}(g)dg.
\end{gather*}
In particular for any $f \in {C}(K \backslash G/K)$,  
\begin{align*}
    [I_{R_{\ell}}f](t) = \int_{T}\!R_{\ell}(gK, t)f(gK)\mu(dgK) = \int_{G}\mathfrak{f}(g)\varphi_{\ell}(g^{-1}g_{t})dg = \mathfrak{f} * \varphi_{\ell}(g_{t}).
\end{align*}
Therefore, since $(G,K)$ is Gelfand,
\begin{align*}
    [I_{R_{2}}(I_{R_{1}}f)](t) &= (\mathfrak{f} * \varphi_{1}) * \varphi_{2}(g_{t}) \\
    &= \mathfrak{f} * (\varphi_{1} * \varphi_{2})(g_{t}) \\
    &= (\mathfrak{f} * \varphi_{2}) * \varphi_{1}(g_{t}) = [I_{R_{1}}(I_{R_{2}}f)](t).
\end{align*}
That is, the integral operators $I_{R_{1}}$ and $I_{R_{2}}$ commute. In other words, the symmetry of $T$, in alignment with $R_{\ell}$, determines the symmetry of $I_{R_{\ell}}$. For the remaining, we want to discuss how Mercer's theorem coupled with the theory of commuting integral operators allows us to describe equivalent Gaussian measures for potentially less symmetric spaces. In particular, $X$ does not necessarily need to be stationary under $P_{\ell}$. Therefore, let $T$ and $R_{\ell}$ be again as at the beginning of this sections. A key ingredient is the following theorem, which is the classical Feldman--H\'{a}jek theorem (see for instance Bogachev's~\cite{Bogachev}), when viewed through the lens of Mercer's theorem. Given the intuition developed in Example~\ref{example:CovSphere}, we present a proof using an RKHS approach.
\begin{theorem}[Feldman--H\'{a}jek--Mercer]\label{thm:MercerCommIntegralOP}
    Assume that $I_{R_{1}}$ and $I_{R_{2}}$ commute. Then $P_{1} \equiv P_{2}$ on $\sigma_{T}(X)$ if and only if $\sum_{n}\big(1-\nicefrac{\lambda_{n}(1)}{\lambda_{n}(2)}\big)^{2} < \infty$, where $\{\lambda_{n}(\ell)\}$ are the eigenvalues of $I_{R_{\ell}}$ for a common set of eigenfunctions $\{e_{n}\}$ as given by Mercer's theorem.   
\end{theorem}

\begin{remark}\label{remark:SymDiagCommute}
    Since $I_{R_{\ell}}$ is compact and self-adjoint, the assumption that $I_{R_{1}}$ and $I_{R_{2}}$ commute is equivalent to their simultaneous diagonalizability. In particular, the assumption of Theorem~\ref{thm:MercerCommIntegralOP} can be substituted with the latter equivalent statement.  
\end{remark}

\begin{proof}[Proof of Theorem~\ref{thm:MercerCommIntegralOP}]
    By the previous remark, we deduce from Mercer's theorem that
    \begin{equation} \label{eq:KernelInMercerBasis}
        R_{\ell}(s,t) = \sum_{n}\lambda_{n}(\ell)e_{n}(s)e_{n}(t),
    \end{equation}
    with $\{\lambda_{n}(\ell)\}$ and $\{e_{n}\}$ as in the statement of the theorem. In particular, $\{\sqrt{\lambda_{n}(\ell)}e_{n}\}$ is an orthonormal basis of $H_{T}(R_{\ell})$. From this point forward, we apply the reasoning developed in Example~\ref{example:CovSphere}. If $P_{1} \equiv P_{2}$ on $\sigma_{T}(X)$, it follows from Theorem~\ref{thm:Aronszajn-Neveu} that $\mathcal{D}$ belongs to $H_{T}(R_{1}^{2\otimes})$ and $H_{T}(R_{1}) = H_{T}(R_{2})$. We express $\mathcal{D}(s,t) = \sum_{i,j}\alpha_{i,j}\sqrt{\lambda_{n}(2)\lambda_{n}(2)}e_{i}(s)e_{j}(t)$, where $\lVert \mathcal{D} \rVert_{R_{2}^{\otimes}}^{2} = \sum_{i,j}\alpha_{i,j}^{2}$. Comparing this expansion with \eqref{eq:KernelInMercerBasis} identifies  
    \begin{equation*} \label{eq:L2CoefficientsDelta}
        \alpha_{i,j} = \begin{cases} \frac{\lambda_{n}(2) - \lambda_{n}(1)}{\lambda_{n}(2)}, & i = j, \\ 0, & \text{otherwise.}\end{cases}
    \end{equation*}
    In particular, the series 
    \begin{equation}\label{eq:SeriesMercerEigenvalueDiff}
        \sum_{n}\bigg(1-\frac{\lambda_{n}(1)}{\lambda_{n}(2)}\bigg)^{2}
    \end{equation}
    converges. On the other hand, if \eqref{eq:SeriesMercerEigenvalueDiff} converges, $\mathcal{D} \in H_{T}(R_{2}^{2\otimes})$ follows from the definition of the product RKHS. Moreover, the convergence of \eqref{eq:SeriesMercerEigenvalueDiff} implies that the sequences $(\nicefrac{\lambda_{n}(1)}{\lambda_{n}(2)})$ and $(\nicefrac{\lambda_{n}(2)}{\lambda_{n}(1)})$ converge to $1$, so that $H_{T}(R_{1}) = H_{T}(R_{2})$ and therefore, $\mathcal{D} \in H_{T}(R_{1}^{2\otimes})$. Applying Theorem~\ref{thm:Aronszajn-Neveu}, this completes the argument.
\end{proof}
If the interplay between $T$ and $R_{\ell}$ does not guarantee the commutativity of $I_{R_{1}}$ and $I_{R_{2}}$, we can still enforce it through a suitable choice of the covariance model. We adopt the approach proposed by Narcowich~\cite{Narcowich}, who investigates strictly positive-definite kernels on compact Riemannian manifolds without boundary that expand by means of an orthonormal basis of eigenfunctions of the Laplace--Beltrami operator. For an introduction to the terminology, see, for instance, the books by Lee~\cite{Lee} or Sakai~\cite{Sakai}. Let $T=M$ be a compact, connected, orientable, $m$-dimensional $C^{\infty}$ Riemannian manifold without boundary and $C^{\infty}$ metric $g_{ij}$. As usual, we denote by $g^{ij}$ the components of the inverse of $g_{ij}$ and write $g = \operatorname{det}g_{ij}$. The measure $\mu$ is taken to be the Riemannian volume form on $M$ and $C^{\infty}(M)$ denotes the space of smooth functions defined on $M$. Examples include K3 surfaces, surfaces of genus $g$, the Grassmannian $\operatorname{Gr}(2,4)$ and toric manifolds. In local coordinates, the Laplace--Beltrami operator $\Delta \colon C^{\infty}(M) \to C^{\infty}(M)$ takes the form
\begin{gather*}
    \Delta f = g^{-1/2}\sum_{i,j=1}^{m}\frac{\partial}{\partial x^{i}}\bigg(g^{1/2}g^{ij}\frac{\partial f}{\partial x^{j}}\bigg).
\end{gather*}
It is known that the eigenvalues of $\Delta$ form an increasing and unbounded sequence of nonnegative real numbers $\lambda_{1} \leq \lambda_{2} \leq \cdots$. Moreover, the eigenspaces are finite-dimensional and their direct product is equal to $L^{2}(\mu)$. In particular, it can be assumed w.l.o.g.\ that the set of eigenfunctions $\{F_{n}\}$ constitutes an orthonormal basis of $L^{2}(\mu)$. 

Narcowich~\cite{Narcowich} studies functions $K \colon M \times M \to \mathbb{R}$, $K \in L^{2}(\mu\otimes\mu)$, of the following type:
\begin{equation}\label{eq:ManifoldKernel}
    K(s,t) = \sum_{n}a_{n}F_{n}(s)F_{n}(t),\quad \sum_{n} a_{n}^{2}< \infty.
\end{equation}
Actually, he puts emphasis on smooth kernels (cf.\ his Proposition~3.5 and Theorem~3.6). In the context of Mercer kernels, Theorem~\ref{thm:MercerCommIntegralOP} straightforwardly characterizes pairs of equivalent Gaussian distributions whose covariance functions expand according to \eqref{eq:ManifoldKernel}.

\begin{corollary}\label{corollary:EquivOrthoGaussCompactManifolds}
    Let $R_{1}(s,t) = \sum_{n}a_{n}(1)F_{n}(s)F_{n}(t)$ and $R_{2}(s,t) = \sum_{n}a_{n}(2)F_{n}(s)F_{n}(t)$ be of type \eqref{eq:ManifoldKernel}. Then $P_{1} \equiv P_{2}$ on $\sigma_{M}(X)$ if and only if $\sum_{n}\big(1-\nicefrac{a_{n}(1)}{a_{n}(2)}\big)^{2} < \infty$. 
\end{corollary}

\begin{proof}
    By construction, $\{F_{n}\}$ is an orthonormal basis of $L^{2}(\mu)$. Hence, the integral operator associated with $R_{\ell}$ has eigenfunctions $\{F_{n}\}$ with corresponding eigenvalues $\{a_{n}(\ell)\}$. Since $I_{R_{1}}$ and $I_{R_{2}}$ commute (cf.\ Remark~\ref{remark:SymDiagCommute}), we apply Theorem~\ref{thm:MercerCommIntegralOP} and the proof of the corollary is complete.
\end{proof}

Viewed from a classical Feldman--H\'{a}jek perspective, Corollary~\ref{corollary:EquivOrthoGaussCompactManifolds} has already been applied to specific choices of covariance models, including squared exponential and Mat\'{e}rn kernels (cf.\ Li et al.~\cite{MaternCompactRiemannianManifold}, following the work of Borovitskiy et al.~\cite{NEURIPS2020}, who provide $\{a_{n}(\ell)\}$ for the latter two correlation models\footnote{In the standard Euclidean setting, where $T$ is a bounded subset of $\Rd$ ($d=1,2,3$), the corresponding analysis is done by Zhang~\cite{zhang}. In the special case of the exponential kernel, for $T = [0,b]$, a compact interval on the line, the result goes back to Striebel \cite{Striebel} (cf.\ \cite{Capon}).}). As a simple illustrative example, consider the heat kernel. 

\begin{example}
    Let $\kappa_{\ell}$ be strictly positive and define
    \begin{gather*}
        K_{\ell}(s,t) = \sum_{n}a_{n}(\ell)F_{n}(s)F_{n}(t), \quad a_{n}(\ell) = \operatorname{e}^{-\lambda_{n}\kappa_{\ell}}, \quad n \in \mathbb{N}.
    \end{gather*}
    That is $K_{\ell}$ takes the form of a heat kernel (see for instance Grigor’yan~\cite{Grigoryan}). In particular, $K_{\ell}$ is smooth and strictly positive definite. Let $R_{\ell} = K_{\ell}$. We derive from Corollary~\ref{corollary:EquivOrthoGaussCompactManifolds} that for $P_{1} \equiv P_{2}$ on $\sigma_{M}(X)$, it must hold that $\lim_{n \to \infty}\operatorname{e}^{-\lambda_{n}(\kappa_{1}-\kappa_{2})} = 1$. The latter implies that $\kappa_{1} = \kappa_{2}$. In particular, $P_{1} \equiv P_{2}$ on $\sigma_{M}(X)$ if and only if $R_{1} = R_{2}$. Especially, given the present setting, the result does not depend on the choice of $M$.
\end{example}

\section{Remarks on covariance parameter estimation} \label{se:ApplicationGeoStat}

Let $T_{\kern-0.6pt\infty} = \{t_{1}, t_{2}, \dotsc\}$ be a countably infinite subset of $T$. Make the assumption that $P_{1} \perp P_{2}$ on $\sigma_{T}(X)$ if and only if $P_{1} \perp P_{2}$ on $\sigma_{T_{\kern-0.6pt\infty}}(X)$ (e.g.\ if $T_{\kern-0.6pt\infty}$ is dense in $T$ and $X$ has continuous sample paths on $T$). Let $\Theta \subset \mathbb{R}^{p}$ and take $P_{\theta}$, $\theta \in \Theta$, to be a parametric family of centered Gaussian distributions on $\sigma_{T_{\kern-0.6pt\infty}}(X)$ with strictly positive-definite covariance function $R_{\theta}$. We can write $\sigma_{T_{\kern-0.6pt\infty}}(X) = \sigma\big(\cup_{n}\sigma(Y_{n})\big)$, with $\sigma(Y_{n})$ the $\sigma$-filed generated by the random vector $Y_{n} = (X_{t_{1}}, \dotsc, X_{t_{n}})$. The set $\Theta$ is regarded as the parameter space and any sequence of random variables $(\hat{\theta}_{n})$ which maximizes the likelihood function over the parameter space is referred to as a sequence of maximum likelihood (ML) estimators. In the present case, the likelihood function is given by $\theta \mapsto p_{\theta}^{n}(Y_{n})$, where $p_{\theta}^{n}$ is the $n$-dimensional normal density of the distribution $P_{\theta}$ restricted to $\sigma(Y_{n})$. Given $\theta_{0} \in \Theta$, a sequence of ML estimators is said to be strongly consistent for $\theta_{0}$ if 
\begin{gather*}
    P_{0}\big(\hat{\theta}_{n} \xrightarrow[]{ n \to \infty } \theta_{0}\big) = 1, \quad P_{0} = P_{\theta_{0}}.
\end{gather*}
It is said to be weakly consistent with limit $\theta_{0}$, if $(\hat{\theta}_{n})$ converges to $\theta_{0}$ in probability $P_{0}$. It turns out that the separability condition
\begin{equation} \label{eq:SeparabilityCondition}
    \theta_{1} \neq \theta_{2} \ \Rightarrow \ P_{\theta_{1}} \perp P_{\theta_{2}} \text{ on $\sigma_{T_{\kern-0.6pt\infty}}(X)$,} \quad \theta_{1},\theta_{2} \in \Theta,
\end{equation}
relates to the feasibility to estimate $\theta_{0}$ consistently. In particular, a violation of \eqref{eq:SeparabilityCondition} can lead to inconsistent ML estimators \cite{zhang,zhang2005towards, anderes2010consistent, bevilacqua2019estimation,Bachoc2021}. Building on a classical likelihood argument (cf.\ Wald~\cite{wald}), it is possible to obtain strongly consistent ML covariance parameter estimators for families of orthogonal Gaussian distributions. The argument relies on the fact that under the separability condition \eqref{eq:SeparabilityCondition}, the Radon--Nikodym derivative $\nicefrac{p_{\theta}^{n}(Y_{n})}{p_{0}^{n}(Y_{n})}$ converges to zero with $P_{0}$ probability one whenever $\theta \neq \theta_{0}$ (cf.\ Theorem~1 on p.\ 442 in \cite{StochasticProcesses}). Then, under suitable conditions on the parameter space and the likelihood function, the latter convergence can be strengthened to a $P_{0}$ a.s.\ uniform convergence on $\Theta$ outside a small neighborhood of $\theta_{0}$ (see for instance Zhang~\cite{zhang}). Regarding weak consistency, the reader is referred to \cite{sweeting1980uniform,mardia1984maximum,bachoc2014asymptotic}. In particular, in the latter two, weakly consistent estimators are obtained under the assumption that the distance between coordinates from $T$ is uniformly bounded away from zero.

%% file: references.bib
@PREAMBLE{ " \newcommand{\noopsort}{} " }

@book{GaussianRandomProcesses,
   AUTHOR = {Ibragimov, Ildar A. and Rozanov, Yuri A.},
    TITLE = {Gaussian Random Processes},
PUBLISHER = {Springer, New York},
     YEAR = {1978},
}

@book {RKHSBook,
    AUTHOR = {Paulsen, Vern I. and Raghupathi, Mrinal},
     TITLE = {An introduction to the theory of reproducing kernel {H}ilbert
              spaces},
    SERIES = {Cambridge Studies in Advanced Mathematics},
    VOLUME = {152},
 PUBLISHER = {Cambridge University Press, Cambridge},
      YEAR = {2016},
     PAGES = {x+182},
  MRNUMBER = {3526117},
}

@book{Jeffreys,
    AUTHOR = {Jeffreys, Harold},
     TITLE = {Theory of probability},
   EDITION = {Second},
 PUBLISHER = {Oxford},
      YEAR = {1948},
}

@book{Wendland_Book,
    AUTHOR = {Wendland, Holger},
     TITLE = {Scattered data approximation},
    SERIES = {Cambridge Monographs on Applied and Computational Mathematics},
    VOLUME = {17},
 PUBLISHER = {Cambridge University Press, Cambridge},
      YEAR = {2005},
     PAGES = {x+336},
  MRNUMBER = {2131724},
}

@book{Bogachev,
    AUTHOR = {Bogachev, Vladimir I.},
     TITLE = {Gaussian measures},
    SERIES = {Mathematical Surveys and Monographs},
    VOLUME = {62},
 PUBLISHER = {American Mathematical Society, Providence, RI},
      YEAR = {1998},
     PAGES = {xii+433},
  MRNUMBER = {1642391},
}

@book{Neveu,
    AUTHOR = {Neveu, Jacques},
     TITLE = {Processus al\'{e}atoires {G}aussiens},
 PUBLISHER = {Les Presses de l'Universit\'{e} de Montr\'{e}al, Canada},
      YEAR = {1968},
}

@book{Yadrenko,
    AUTHOR = {Yadrenko, Mikhail I.},
     TITLE = {Spectral Theory of Random Fields},
 PUBLISHER = {Optimization Software, Inc., New York},
      YEAR = {1983},
}

@book{Erdelyi,
    AUTHOR = {Erd\'{e}lyi, Arthur and Magnus, Wilhelm and Oberhettinger, Fritz
              and Tricomi, Francesco G.},
     TITLE = {Higher Transcendental Functions, {V}ol. {II}},
      NOTE = {Based on notes left by Harry Bateman,
              Reprint of the 1953 original},
 PUBLISHER = {Robert E. Krieger Publishing Co., Inc., Melbourne, Fla.},
      YEAR = {1981},
}

@book{Malyarenko,
    AUTHOR = {Malyarenko, Anatoliy},
     TITLE = {Invariant random fields on spaces with a group action},
    SERIES = {Probability and its Applications},
      NOTE = {With a foreword by Nikolai Leonenko},
 PUBLISHER = {Springer, Heidelberg},
      YEAR = {2013},
     PAGES = {xviii+261},
  MRNUMBER = {2977490},
}

@article {Aronszajn,
    AUTHOR = {Aronszajn, Nachman},
     TITLE = {Theory of reproducing kernels},
   JOURNAL = {Trans. Amer. Math. Soc.},
  FJOURNAL = {Transactions of the American Mathematical Society},
    VOLUME = {68},
      YEAR = {1950},
     PAGES = {337--404},
  MRNUMBER = {51437},
}

@inproceedings{KallianpurOodaira1,
    AUTHOR = {Kallianpur, Gopinath and Oodaira, Hiroshi},
     TITLE = {The equivalence and singularity of {G}aussian measures},
 BOOKTITLE = {Proc. {S}ympos. {T}ime {S}eries {A}nalysis},
     PAGES = {279--291},
 PUBLISHER = {Wiley, New York},
      YEAR = {1963},
  MRNUMBER = {149527},
}

@book{Oodaira,
    AUTHOR = {Oodaira, Hiroshi},
     TITLE = {The equivalence of {G}aussian stochastic processes},
      NOTE = {Thesis for the degree of Ph.D.},
 PUBLISHER = {Michigan State University},
      YEAR = {1963},
  MRNUMBER = {2613909},
}

@book {Wolf,
    AUTHOR = {Wolf, Joseph A.},
     TITLE = {Harmonic analysis on commutative spaces},
    SERIES = {Mathematical Surveys and Monographs},
    VOLUME = {142},
 PUBLISHER = {American Mathematical Society},
      YEAR = {2007},
     PAGES = {xvi+387},
  MRNUMBER = {2328043},
}

@article {HajekJDiv,
    AUTHOR = {H\'{a}jek, Jaroslav},
     TITLE = {A property of {$J$}-divergences of marginal probability
              distributions},
   JOURNAL = {Czechoslovak Math. J.},
  FJOURNAL = {Czechoslovak Mathematical Journal},
    VOLUME = {8},
      YEAR = {{\vphantom{a}}1958a},
     PAGES = {460--463},
  MRNUMBER = {99712},
}

@inproceedings {ParzenRKHS,
    AUTHOR = {Parzen, Emanuel},
     TITLE = {Probability density functionals and reproducing kernel
              {H}ilbert spaces},
 BOOKTITLE = {{T}ime {S}eries {A}nalysis},
     PAGES = {155--169},
 PUBLISHER = {Wiley, New York-London},
      YEAR = {1963},
  MRNUMBER = {149634},
}

@article {Parzen1962,
    AUTHOR = {Parzen, Emanuel},
     TITLE = {Extraction and detection problems and reproducing kernel
              {H}ilbert spaces},
   JOURNAL = {J. SIAM Control Ser. A},
  FJOURNAL = {Journal of the Society for Industrial and Applied Mathematics.
              Series A. On Control},
    VOLUME = {1},
      YEAR = {1962},
     PAGES = {35--62},
}

@techreport{ParzenInferenceTimeSeriesRKHSMethods,
  author = {Parzen, Emanuel},
  title = {Statistical {I}nference on {T}ime {S}eries by {H}ilbert {S}pace {M}ethods, {I}},
  institution = {Applied Mathematics and Statistics Laboratory, Stanford University},
  year = {1959},
  number = {No. 23}
}

@inproceedings {ParzenTimeSeries,
    AUTHOR = {Parzen, Emanuel},
     TITLE = {Regression analysis of continuous parameter time series},
 BOOKTITLE = {{P}roceedings of the {F}ourth {B}erkeley {S}ymposium on {M}athematical {S}tatistics and {P}robability,
              {V}olume 1: {C}ontributions to the {T}heory of {S}tatistics},
     PAGES = {469--489},
    VOLUME = {4},
 PUBLISHER = {University of California Press},
      YEAR = {1961},
  MRNUMBER = {146931},
}

@book{SteinWeiss,
    AUTHOR = {Stein, Elias M. and Weiss, Guido},
     TITLE = {Introduction to {F}ourier analysis on {E}uclidean spaces},
    SERIES = {Princeton Mathematical Series},
 PUBLISHER = {Princeton University Press},
      YEAR = {1971},
     PAGES = {x+297},
  MRNUMBER = {304972},
}

@book{Rudin,
    AUTHOR = {Rudin, Walter},
     TITLE = {Fourier analysis on groups},
    SERIES = {Interscience Tracts in Pure and Applied Mathematics, No. 12},
 PUBLISHER = {Interscience Publishers (a division of John Wiley \& Sons), New York-London},
      YEAR = {1962},
     PAGES = {ix+285},
  MRNUMBER = {152834},
}

@article{HajekOriginal,
    AUTHOR = {H\'{a}jek, Jaroslav},
     TITLE = {On a property of normal distributions of any stochastic
              process},
   JOURNAL = {Czechoslovak Math. J.},
  FJOURNAL = {Czechoslovak Mathematical Journal},
    VOLUME = {8},
      YEAR = {{\vphantom{b}}1958b},
     PAGES = {610--618},
  MRNUMBER = {104290},
}

@article {Capon,
    AUTHOR = {Capon, Jack},
     TITLE = {Radon--{N}ikodym derivatives of stationary {G}aussian measures},
   JOURNAL = {Ann. Math. Statist.},
  FJOURNAL = {Annals of Mathematical Statistics},
    VOLUME = {35},
      YEAR = {1964},
     PAGES = {517--531},
  MRNUMBER = {161377},
}

@incollection {Chatterji_Mandrekar,
    AUTHOR = {Chatterji, Srishti D. and Mandrekar, Vidyadhar},
     TITLE = {Equivalence and singularity of {G}aussian measures and
              applications},
 BOOKTITLE = {Probabilistic analysis and related topics, {V}ol. 1},
     PAGES = {169--197},
 PUBLISHER = {Academic Press, New York-London},
      YEAR = {1978},
  MRNUMBER = {478320},
}

@article{Rozanov,
author = {Rozanov, Yuri A.},
title = {On the density of one {G}aussian measure with respect to another},
journal = {Teor. Veroyatnost. i Primenen.},
volume = {7},
pages = {84--89 (in Russian)},
year = {1962},
note = {English transl.: Theory Probab. Appl. {\bf 7}, (1962), 82--87}
}

@article {Feldman1958,
    AUTHOR = {Feldman, Jacob},
     TITLE = {Equivalence and perpendicularity of {G}aussian processes},
   JOURNAL = {Pacific J. Math.},
  FJOURNAL = {Pacific Journal of Mathematics},
    VOLUME = {8},
      YEAR = {1958},
     PAGES = {699--708},
  MRNUMBER = {102760},
}

@incollection{Bachoc2021,
   AUTHOR    = {Bachoc, Fran\c{c}ois},
   EDITOR    = {Daouia, Abdelaati and Ruiz-Gazen, Anne},
   TITLE     = {Asymptotic Analysis of Maximum Likelihood Estimation of Covariance Parameters for {G}aussian Processes: {A}n Introduction with Proofs},
   BOOKTITLE = {Advances in Contemporary Statistics and Econometrics},
   YEAR      = {2021},
   PUBLISHER = {Springer},
   ADDRESS   = {Cham},
   PAGES     = {283--303},
   MRNUMBER = {4299285},
}

@incollection{Godement1958,
     author = {Godement, Roger},
     title = {Introduction aux travaux de {A.} {Selberg}},
     booktitle = {S\'eminaire Bourbaki : ann\'ees 1956/57 - 1957/58, expos\'es 137-168},
     series = {Ast\'erisque},
     note = {talk:144},
     pages = {95--110},
     publisher = {Soci\'et\'e math\'ematique de France},
     number = {4},
     year = {1958},
     mrnumber = {1610957},
}

@article {GodementSpherical,
    AUTHOR = {Godement, Roger},
     TITLE = {A theory of spherical functions. {I}},
   JOURNAL = {Trans. Amer. Math. Soc.},
  FJOURNAL = {Transactions of the American Mathematical Society},
    VOLUME = {73},
      YEAR = {1952},
     PAGES = {496--556},
  MRNUMBER = {52444},
}

@article {wald,
    AUTHOR = {Wald, Abraham},
     TITLE = {Note on the consistency of the maximum likelihood estimate},
   JOURNAL = {Ann. Math. Statistics},
    VOLUME = {20},
      YEAR = {1949},
     PAGES = {595--601},
  MRNUMBER = {32169},
}

@article {zhang,
    AUTHOR = {Zhang, Hao},
     TITLE = {Inconsistent estimation and asymptotically equal
              interpolations in model-based geostatistics},
   JOURNAL = {J. Amer. Statist. Assoc.},
    VOLUME = {99},
      YEAR = {2004},
     PAGES = {250--261},
  MRNUMBER = {2054303},
}

@book {StochasticProcesses,
    AUTHOR = {Gikhman, Iosif I. and Skorokhod, Anatoli V.},
     TITLE = {The theory of stochastic processes {I}},
    SERIES = {Classics in Mathematics},
 PUBLISHER = {Springer-Verlag},
   ADDRESS = {Berlin},
      YEAR = {2004},
     PAGES = {viii+574},
}

@article {anderes2010consistent,
    AUTHOR = {Anderes, Ethan},
     TITLE = {On the consistent separation of scale and variance for
              {G}aussian random fields},
   JOURNAL = {Ann. Math. Statistics},
  FJOURNAL = {The Annals of Statistics},
    VOLUME = {38},
      YEAR = {2010},
    NUMBER = {2},
     PAGES = {870--893},
  MRNUMBER = {2604700},
}

@article {bevilacqua2019estimation,
    AUTHOR = {Bevilacqua, Moreno and Faouzi, Tarik and Furrer, Reinhard and
              Porcu, Emilio},
     TITLE = {Estimation and prediction using generalized {W}endland
              covariance functions under fixed domain asymptotics},
   JOURNAL = {Ann. Math. Statistics},
  FJOURNAL = {The Annals of Statistics},
    VOLUME = {47},
      YEAR = {2019},
    NUMBER = {2},
     PAGES = {828--856},
  MRNUMBER = {3909952},
}

@article {zhang2005towards,
    AUTHOR = {Zhang, Hao and Zimmerman, Dale L.},
     TITLE = {Towards reconciling two asymptotic frameworks in spatial
              statistics},
   JOURNAL = {Biometrika},
  FJOURNAL = {Biometrika},
    VOLUME = {92},
      YEAR = {2005},
    NUMBER = {4},
     PAGES = {921--936},
  MRNUMBER = {2234195},
}

@article {GaussianSphereEquivOrtho,
    AUTHOR = {Arafat, Ahmed and Porcu, Emilio and Bevilacqua, Moreno and
              Mateu, Jorge},
     TITLE = {Equivalence and orthogonality of {G}aussian measures on
              spheres},
   JOURNAL = {J. Multivariate Anal.},
  FJOURNAL = {Journal of Multivariate Analysis},
    VOLUME = {167},
      YEAR = {2018},
     PAGES = {306--318},
  MRNUMBER = {3830648},
}

@article {Schoenberg,
    AUTHOR = {Schoenberg, Isaac J.},
     TITLE = {Positive definite functions on spheres},
   JOURNAL = {Duke Math. J.},
  FJOURNAL = {Duke Mathematical Journal},
    VOLUME = {9},
      YEAR = {1942},
     PAGES = {96--108},
  MRNUMBER = {5922},
}

@article{CameronMartin,
    AUTHOR = {Cameron, Robert H. and Martin, William T.},
     TITLE = {The behavior of measure and measurability under change of
              scale in {W}iener space},
   JOURNAL = {Bull. Amer. Math. Soc.},
  FJOURNAL = {Bulletin of the American Mathematical Society},
    VOLUME = {53},
      YEAR = {1947},
     PAGES = {130--137},
  MRNUMBER = {19259},
}

@article{GneitingMultiquadric,
    AUTHOR = {Gneiting, Tilmann},
     TITLE = {Strictly and non-strictly positive definite functions on
              spheres},
   JOURNAL = {Bernoulli},
  FJOURNAL = {Bernoulli. Official Journal of the Bernoulli Society for
              Mathematical Statistics and Probability},
    VOLUME = {19},
      YEAR = {2013},
    NUMBER = {4},
     PAGES = {1327--1349},
  MRNUMBER = {3102554},
}

@article {Tamagawa,
    AUTHOR = {Tamagawa, Tsuneo},
     TITLE = {On {S}elberg's trace formula},
   JOURNAL = {J. Fac. Sci. Univ. Tokyo Sect. I},
  FJOURNAL = {Journal of the Faculty of Science. University of Tokyo.
              Section I},
    VOLUME = {8},
      YEAR = {1960},
     PAGES = {363--386},
  MRNUMBER = {123633},
}

@article {Striebel,
    AUTHOR = {Striebel, Charlotte T.},
     TITLE = {Densities for stochastic processes},
   JOURNAL = {Ann. Math. Statist.},
  FJOURNAL = {Annals of Mathematical Statistics},
    VOLUME = {30},
      YEAR = {1959},
     PAGES = {559--567},
  MRNUMBER = {104330},
}

@article {mardia1984maximum,
    AUTHOR = {Mardia, Kanti V. and Marshall, Roger J.},
     TITLE = {Maximum likelihood estimation of models for residual
              covariance in spatial regression},
   JOURNAL = {Biometrika},
  FJOURNAL = {Biometrika},
    VOLUME = {71},
      YEAR = {1984},
    NUMBER = {1},
     PAGES = {135--146},
  MRNUMBER = {738334},
}

@article {sweeting1980uniform,
    AUTHOR = {Sweeting, Trevor J.},
     TITLE = {Uniform asymptotic normality of the maximum likelihood
              estimator},
   JOURNAL = {The Annals of Statistics},
  FJOURNAL = {The Annals of Statistics},
    VOLUME = {8},
      YEAR = {1980},
    NUMBER = {6},
     PAGES = {1375--1381},
  MRNUMBER = {594652},
}

@article{bachoc2014asymptotic,
    AUTHOR = {Bachoc, Fran\c{c}ois},
     TITLE = {Asymptotic analysis of the role of spatial sampling for
              covariance parameter estimation of {G}aussian processes},
   JOURNAL = {Journal of Multivariate Analysis},
  FJOURNAL = {Journal of Multivariate Analysis},
    VOLUME = {125},
      YEAR = {2014},
     PAGES = {1--35},
  MRNUMBER = {3163828},
}

@article {KullbackLeibler,
    AUTHOR = {Kullback, Solomon and Leibler, Richard A.},
     TITLE = {On information and sufficiency},
   JOURNAL = {Ann. Math. Statistics},
  FJOURNAL = {Annals of Mathematical Statistics},
    VOLUME = {22},
      YEAR = {1951},
     PAGES = {79--86},
  MRNUMBER = {39968},
}

@article {Bochner,
    AUTHOR = {Bochner, Salomon},
     TITLE = {Monotone {F}unktionen, {S}tieltjessche {I}ntegrale und
              harmonische {A}nalyse},
   JOURNAL = {Math. Ann.},
  FJOURNAL = {Mathematische Annalen},
    VOLUME = {108},
      YEAR = {1933},
    NUMBER = {1},
     PAGES = {378--410},
  MRNUMBER = {1512856},
}

@book{Helgason,
    AUTHOR = {Helgason, Sigurdur},
     TITLE = {Groups and geometric analysis},
    SERIES = {Mathematical Surveys and Monographs},
    VOLUME = {83},
 PUBLISHER = {American Mathematical Society, Providence, RI},
      YEAR = {2000},
     PAGES = {xxii+667},
  MRNUMBER = {1790156},
}

@article {Gelfand,
    AUTHOR = {Gelfand, Israel M.},
     TITLE = {Spherical functions in symmetric {R}iemann spaces},
   JOURNAL = {Doklady Akad. Nauk. SSSR.},
  FJOURNAL = {Doklady Akad. Nauk. SSSR.},
    VOLUME = {70},
      YEAR = {1950},
     PAGES = {5--8},
  MRNUMBER = {33832},
}

@article {Cartan,
    AUTHOR = {Cartan, \'{E}lie},
     TITLE = {Sur la d\'{e}termination d’un syst\`{e}me orthogonal Complet dans un espace de {R}iemann sym\'{e}trique Clos},
   JOURNAL = {Rend. Circ. Mat. Palermo},
  FJOURNAL = {Rendiconti del Circolo Matematico di Palermo},
    VOLUME = {53},
      YEAR = {1929},
     PAGES = {217--252},
}

@article {Weyl,
    AUTHOR = {Weyl, Hermann},
     TITLE = {Harmonics on homogeneous manifolds},
   JOURNAL = {Ann. of Math.},
  FJOURNAL = {Annals of Mathematics. Second Series},
    VOLUME = {35},
      YEAR = {1934},
    NUMBER = {3},
     PAGES = {486--499},
  MRNUMBER = {1503175},
}

@inproceedings {YaglomInvariant,
    AUTHOR = {Yaglom, Akiva M.},
     TITLE = {Second-order homogeneous random fields},
 BOOKTITLE = {Proc. 4th {B}erkeley {S}ympos. {M}ath. {S}tatist. and {P}rob.,
              {V}ol. {II}},
     PAGES = {593--622},
 PUBLISHER = {Univ. California Press, Berkeley-Los Angeles, Calif.},
      YEAR = {1961},
  MRNUMBER = {146880},
}

@article {chow1,
    AUTHOR = {Chow, Thomas R.},
     TITLE = {The equivalence of group invariant positive definite
              functions},
   JOURNAL = {Pacific J. Math.},
  FJOURNAL = {Pacific Journal of Mathematics},
    VOLUME = {26},
      YEAR = {1968},
     PAGES = {25--38},
  MRNUMBER = {231217},
}

@article {chow2,
    AUTHOR = {Chow, Thomas R.},
     TITLE = {Equivalence of {G}aussian stationary processes},
   JOURNAL = {Ann. Math. Statist.},
  FJOURNAL = {Annals of Mathematical Statistics},
    VOLUME = {40},
      YEAR = {1969},
     PAGES = {197--202},
  MRNUMBER = {235610},
}

@article {Grenander,
    AUTHOR = {Grenander, Ulf},
     TITLE = {Stochastic processes and statistical inference},
   JOURNAL = {Ark. Mat.},
  FJOURNAL = {Arkiv f\"{o}r Matematik},
    VOLUME = {1},
      YEAR = {1950},
     PAGES = {195--277},
  MRNUMBER = {39202},
}

@book {Saburou_Yoshihiro_RKHS,
    AUTHOR = {Saitoh, Saburou and Sawano, Yoshihiro},
     TITLE = {Theory of reproducing kernels and applications},
    SERIES = {Developments in Mathematics},
    VOLUME = {44},
 PUBLISHER = {Springer, Singapore},
      YEAR = {2016},
     PAGES = {xviii+452},
  MRNUMBER = {3560890},
}

@article {Askey_Bingham,
    AUTHOR = {Askey, R. and Bingham, Nicholas H.},
     TITLE = {Gaussian processes on compact symmetric spaces},
   JOURNAL = {Z. Wahrscheinlichkeitstheorie und Verw. Gebiete},
  FJOURNAL = {Zeitschrift f\"{u}r Wahrscheinlichkeitstheorie und Verwandte
              Gebiete},
    VOLUME = {37},
      YEAR = {1976/77},
    NUMBER = {2},
     PAGES = {127--143},
  MRNUMBER = {423000},
}

@book{Folland,
    AUTHOR = {Folland, Gerald B.},
     TITLE = {A course in abstract harmonic analysis},
    SERIES = {Textbooks in Mathematics},
   EDITION = {Second},
 PUBLISHER = {CRC Press, Boca Raton, FL},
      YEAR = {2016},
     PAGES = {xiii+305 pp.+loose errata},
  MRNUMBER = {3444405},
}

@article {Narcowich,
    AUTHOR = {Narcowich, Francis J.},
     TITLE = {Generalized {H}ermite interpolation and positive definite
              kernels on a {R}iemannian manifold},
   JOURNAL = {J. Math. Anal. Appl.},
  FJOURNAL = {Journal of Mathematical Analysis and Applications},
    VOLUME = {190},
      YEAR = {1995},
    NUMBER = {1},
     PAGES = {165--193},
  MRNUMBER = {1314111},
}

@book{Lee,
    AUTHOR = {Lee, John M.},
     TITLE = {Introduction to {R}iemannian manifolds},
    SERIES = {Graduate Texts in Mathematics},
    VOLUME = {176},
      NOTE = {Second edition},
 PUBLISHER = {Springer, Cham},
      YEAR = {2018},
     PAGES = {xiii+437},
  MRNUMBER = {3887684},
}

@book {Sakai,
    AUTHOR = {Sakai, Takashi},
     TITLE = {Riemannian geometry},
    SERIES = {Translations of Mathematical Monographs},
    VOLUME = {149},
 PUBLISHER = {American Mathematical Society, Providence, RI},
      YEAR = {1996},
     PAGES = {xiv+358},
  MRNUMBER = {1390760},
}

@inproceedings{NEURIPS2020,
 author = {Borovitskiy, Viacheslav and Terenin, Alexander and Mostowsky, Peter and Deisenroth, Marc P.},
 booktitle = {Advances in Neural Information Processing Systems},
 pages = {12426--12437},
 title = {{M}at\'{e}rn {G}aussian processes on {R}iemannian manifolds},
 volume = {33},
 year = {2020}
}

@book {Grigoryan,
    AUTHOR = {Grigor'yan, Alexander},
     TITLE = {Heat kernel and analysis on manifolds},
    SERIES = {AMS/IP Studies in Advanced Mathematics},
    VOLUME = {47},
 PUBLISHER = {American Mathematical Society, Providence, RI; International
              Press, Boston, MA},
      YEAR = {2009},
     PAGES = {xviii+482},
  MRNUMBER = {2569498},
}

@article{Mercer,
    author = {Mercer, James},
    title = {Functions of positive and negative type, and their connection with the theory of integral equations},
    journal = {Philosophical Transactions of the Royal Society A},
    volume = {209},
    pages = {415--446},
    year = {1909},
    month = {01},
    doi = {10.1098/rsta.1909.0016},
}

@article{MaternCompactRiemannianManifold,
    AUTHOR = {Li, Didong and Tang, Wenpin and Banerjee, Sudipto},
     TITLE = {Inference for {G}aussian processes with {M}at\'{e}rn
              covariogram on compact {R}iemannian manifolds},
   JOURNAL = {J. Mach. Learn. Res.},
  FJOURNAL = {Journal of Machine Learning Research (JMLR)},
    VOLUME = {24},
      YEAR = {2023},
    Number = {101},
     PAGES = {1--26},
  MRNUMBER = {4583262},
}

@article {Selberg,
    AUTHOR = {Selberg, Atle},
     TITLE = {Harmonic analysis and discontinuous groups in weakly symmetric
              {R}iemannian spaces with applications to {D}irichlet series},
   JOURNAL = {J. Indian Math. Soc. (N.S.)},
  FJOURNAL = {The Journal of the Indian Mathematical Society. New Series},
    VOLUME = {20},
      YEAR = {1956},
     PAGES = {47--87},
  MRNUMBER = {88511},
}
